\documentclass[12pt,a4paper]{amsart}
\pdfoutput=1

\usepackage{setspace}
\usepackage[centertags]{amsmath}
\usepackage{amsfonts,ifthen,latexsym,amssymb,amsthm, amscd}

\usepackage[T1]{fontenc}
\usepackage[utf8]{inputenc}
\usepackage[american]{babel}

\usepackage[usenames,dvipsnames]{color}

\definecolor{labelkey}{gray}{.20}
\definecolor{refkey}{gray}{.20}
\definecolor{eqkey}{gray}{.20}

\allowdisplaybreaks

\usepackage[american]{babel} 

\usepackage[perpage]{footmisc}

\usepackage{graphicx}

\usepackage{url}
\usepackage{hyperref}
\hypersetup{colorlinks=true,citecolor=blue,filecolor=blue,linkcolor=blue,urlcolor=blue}

\def\mathclap#1{\text{\hbox to 0pt{\hss$\mathsurround=0pt#1$\hss}}}

\newcounter{results}

\newcounter{questions}

\newtheorem{question}[questions]{Question}

\newcounter{mylabel}

\newcommand{\beq}{\begin{equation}}
\newcommand{\eeq}{\end{equation}}
\newtheorem*{theorem*}{Theorem}
\newtheorem*{lemma*}{Lemma}
\newtheorem*{cory*}{Corollary}
\newtheorem*{question*}{Question}

\newtheorem{theorem}[subsection]{Theorem}

\newtheorem{corollary}[subsection]{Corollary}
\newtheorem{cory}[subsection]{Corollary}
\newtheorem{lemma}[subsection]{Lemma}
\newtheorem{proposition}[subsection]{Proposition}

\newtheorem*{claim*}{Claim}

\theoremstyle{definition}

\theoremstyle{remark}

\newcommand{\al}{{\alpha}}

\newcommand{\be}{{\beta}}

\newcommand{\Om}{{\Omega}}

\newcommand{\eps}{{\varepsilon}}

\newcommand{\De}{{\Delta}}
\newcommand{\ga}{{\gamma}}
\newcommand{\Ga}{{\Gamma}}

\newcommand{\si}{{\sigma}}
\newcommand{\Si}{{\Sigma}}
\renewcommand{\phi}{{\varphi}}

\newcommand\R{\mathbb R}
\newcommand\Q{\mathbb Q}
\newcommand\Z{\mathbb Z}
\newcommand\C{\mathbb C}
\newcommand\N{\mathbb N}

\newcommand{\actson}{\curvearrowright}

\newcommand{\mapsinto}{{\hookrightarrow}}

\newcommand{\wt}[1]{{\widetilde {#1}}}
\newcommand{\wh}[1]{{\widehat {#1}}}
\newcommand{\cal}[1]{{\mathcal #1}}
\newcommand{\call}{\cal}
\makeatletter
\newcommand*\if@single[3]{%
  \setbox0\hbox{${\mathaccent"0362{#1}}^H$}%
  \setbox2\hbox{${\mathaccent"0362{\kern0pt#1}}^H$}%
  \ifdim\ht0=\ht2 #3\else #2\fi
  }
\newcommand*\rel@kern[1]{\kern#1\dimexpr\macc@kerna}
\newcommand*\widebar[1]{\@ifnextchar^{{\wide@bar{#1}{0}}}{\wide@bar{#1}{1}}}
\newcommand*\wide@bar[2]{\if@single{#1}{\wide@bar@{#1}{#2}{1}}{\wide@bar@{#1}{#2}{2}}}
\newcommand*\wide@bar@[3]{%
  \begingroup
  \def\mathaccent##1##2{%
    \if#32 \let\macc@nucleus\first@char \fi
    \setbox\z@\hbox{$\macc@style{\macc@nucleus}_{}$}%
    \setbox\tw@\hbox{$\macc@style{\macc@nucleus}{}_{}$}%
    \dimen@\wd\tw@
    \advance\dimen@-\wd\z@
    \divide\dimen@ 3
    \@tempdima\wd\tw@
    \advance\@tempdima-\scriptspace
    \divide\@tempdima 10
    \advance\dimen@-\@tempdima
    \ifdim\dimen@>\z@ \dimen@0pt\fi
    \rel@kern{0.6}\kern-\dimen@
    \if#31
      \overline{\rel@kern{-0.6}\kern\dimen@\macc@nucleus\rel@kern{0.4}\kern\dimen@}%
      \advance\dimen@0.4\dimexpr\macc@kerna
      \let\final@kern#2%
      \ifdim\dimen@<\z@ \let\final@kern1\fi
      \if\final@kern1 \kern-\dimen@\fi
    \else
      \overline{\rel@kern{-0.6}\kern\dimen@#1}%
    \fi
  }%
  \macc@depth\@ne
  \let\math@bgroup\@empty \let\math@egroup\macc@set@skewchar
  \mathsurround\z@ \frozen@everymath{\mathgroup\macc@group\relax}%
  \macc@set@skewchar\relax
  \let\mathaccentV\macc@nested@a
  \if#31
    \macc@nested@a\relax111{#1}%
  \else
    \def\gobble@till@marker##1\endmarker{}%
    \futurelet\first@char\gobble@till@marker#1\endmarker
    \ifcat\noexpand\first@char A\else
      \def\first@char{}%
    \fi
    \macc@nested@a\relax111{\first@char}%
  \fi
  \endgroup
}
\makeatother

\newcommand*\mcapinn[2]{\vcenter{\hbox{$\mathsurround=0pt
  \ifx\displaystyle#1\textstyle\else#1\fi\bigcap$}}}

\newcommand*\mcupinn[2]{\vcenter{\hbox{$\mathsurround=0pt
  \ifx\displaystyle#1\textstyle\else#1\fi\bigcup$}}}

\newcommand\supp{{\operatorname{supp}}}

\newcommand\tr{\operatorname{tr}}
\newcommand\schr{\operatorname{Schr}}
\newcommand\cay{\operatorname{Cay}}

\newcommand\dimvn{\dim_\text{vN}}

\newcommand\im{\operatorname{im}}

\renewcommand{\ge}{\geqslant}

\newcommand{\Zmod}[1]{\,\Z{/#1\Z}\,}
\newcommand{\bb}[1]{{\mathbb #1}}
\newcommand\st{\operatorname{St}}
\newcommand\Dom{\operatorname{Dom}}
\newcommand{\p}{{\partial}}
\newcommand\Aut{\operatorname{Aut}}

\newcounter{subsecs}[section]
\def\thesubsecs{\alph{subsecs}}
\definecolor{color_toc_subsection}{RGB}{0,0,0}
\newcommand{\subsec}[1]{
    \medskip
    \refstepcounter{subsecs}
    \addcontentsline{toc}{section}{
 	{\color{color_toc_subsection}\small \qquad #1 }
    }
    {\noindent{\scshape{\thesection\thesubsecs}}{ #1}\nopagebreak[4]\smallskip\nopagebreak[4]}
}
\usepackage[parskip=5pt,indention=5pt, scriptsize,tt]{caption}

\begin{document}

\title{On Turing dynamical systems and the Atiyah problem}
\author{{\L}ukasz Grabowski}
\address{University of Warwick, Mathematics Institute, Zeeman Building, Coventry, CV4 7AL, UK}
\begin{abstract}
Main theorems of the article concern the problem of M. Atiyah on possible values of $l^2$-Betti numbers. It is shown that all non-negative real numbers are $l^2$-Betti numbers, and that ``many'' (for example all non-negative algebraic) real numbers are $l^2$-Betti numbers of simply connected manifolds with respect to a free cocompact action. Also an explicit example is constructed which leads to a simply connected  manifold with a transcendental $l^2$-Betti number with respect to an action of the threefold direct product of the lamplighter group $\Zmod{2}\wr \Z$. The main new idea is embedding Turing machines into integral group rings. The main tool developed generalizes known techniques of spectral computations for certain random walk operators to arbitrary operators in groupoid rings of discrete measured groupoids.
\end{abstract}
\maketitle
\addtocontents{toc}{\protect\setstretch{0.3}}
\numberwithin{equation}{section}
\setcounter{tocdepth}{1}
\singlespace
\tableofcontents
\bigskip
\begin{flushright}
 {\textit{To my family}}
\end{flushright}

\section{Introduction}\label{section_intro}
Main theorems of the article concern studying countable discrete groups through so called \textit{$l^2$-Betti numbers}. These are certain invariants, originally introduced by M. Atiyah in \cite{Atiyah1976} to study free cocompact actions of discrete groups on manifolds. Subsequently, they were studied and used in many different contexts in geometry and group theory (e.g. \cite{dodziuk:de_rham_hodge_theory_for_l2_cohomology}, \cite{cheeger_gromov:l2_cohomology_and_group_cohomology}, \cite{gaboriau_invariants_l2_de_relations_dequivalence}). 

A particular question Atiyah asked in \cite{Atiyah1976} was whether $l^2$-Betti numbers can be irrational. Since then, various statements about restrictions on possible values of $l^2$-Betti numbers bear the name \textit{the Atiyah conjecture} (e.g. \cite{dodziuk_linnell_mathai_schick_yates:approximating_l2_invariants}). We depart somewhat from this tradition. Given a countable discrete group $G$, the following question will be referred to as the  \textit{Atiyah problem for $G$}.

\begin{question*}
 What is the set of $l^2$-Betti numbers arising from $G$?
\end{question*}

Let us right away introduce the notation $\call C(G)$ for the set in the above question. Over time it has been realized (see \cite{Eckmann_intro} and \cite{Grigorchuk_Linnel_Schick_Zuk}) that $l^2$-Betti numbers arising from a given group can be defined purely in terms of $G$, without mentioning manifolds. Consequently, the Atiyah problem can also be phrased purely in terms of $G$. This is the approach we adopt in the article and which we now briefly present.

The rational group ring $\Q G$ acts on the Hilbert space $l^2G$ by convolution, and similarly matrices $M_k(\Q G) \cong M_k(\Q)\otimes \Q G$ act on $(l^2G)^k$. We have a trace $\tau_G$ on $\Q G$ defined by $\tau_G(T) : = \langle T\zeta_e, \zeta_e \rangle$, where $\zeta_e\in l^2G$ is the vector corresponding to the neutral element of $G$, and we have the induced trace $\tr\otimes \tau_G$ on $M_k(\Q G)$, also denoted by $\tau_G$. 

Recall that when $R$ is a ${}^*$-ring of operators on a Hilbert space, together with a trace $\tau$ which is normal (i.e., extends in a continuous way to the weak closure of $R$), positive (i.e., $\tau(T^*T)\ge 0$) and faithful (i.e., $\tau(T^*T)=0$ implies $T=0$), then for a self-adjoint $T\in R$ we can compose the usual projection-valued spectral measure with $\tau$, to obtain the (scalar-valued) spectral measure of $T$. In particular, spectral measure of the set $\{0\}$ is called the {\bf von Neumann dimension of the kernel of $T$}, and  is denoted by $\dimvn \ker T$. For a non-self-adjoint $T$, one defines $\dimvn \ker T := \dimvn\ker T^*T$.

It turns out that $\tau_G$ is a positive faithful normal trace on $M_k(\Q G)$ and thus we have a von Neumann dimension. A real number $r$ is said to be an {\bf $l^2$-Betti number arising from $G$} if and only if there exists a matrix $T\in M_k(\Q G)$ such that $\dimvn\ker T = r$. 

Much is known about the Atiyah problem for various particular groups. If $G$ is torsion-free, then $\call C(G)$ is conjectured to be the set of non-negative integers. This statement is known as the \textit{Atiyah conjecture for torsion-free groups} (there is a similar conjecture for groups whose torsion subgroups have bounded orders). Cases for which the Atiyah conjecture is known include elementary amenable groups, free groups (see \cite{Linnel:Atiyah_conjecture_for_free_groups} for both classes)  and braid groups (see \cite{linnell_schick:finite_group_extensions}). Many results follow by applying versions of L{\"u}ck's approximation theorem (see \cite{lueck:aroximating_l2_invariants_by_their_finite}, \cite{dodziuk_linnell_mathai_schick_yates:approximating_l2_invariants}, Lemma 5.2 in \cite{osin_thom_normal_generation_and_l2_betti_numbers_of_groups}) to already established results.  Perhaps the most familiar
consequence of the Atiyah conjecture  is the zero divisors conjecture for torsion-free groups. For other results see \cite{Lueck:Big_book}, Chapter 10.

Before the work of R. Grigorchuk and A. Żuk in \cite{Grigorchuk_Zuk2001}, it had been conjectured that $\call C(G) \subset \Z(\frac1{a_1}, \frac1{a_2},\ldots)$, where $a_i$ are orders of torsion subgroups in $G$. However, in \cite{Grigorchuk_Zuk2001} the authors showed that $\dimvn\ker T = \frac13$ for a certain operator $T$ from the group ring of the lamplighter group $\Zmod{2} \wr \Z$. Recall that the latter group is a semi-direct product of $\Zmod{2}^{\oplus \Z}$ with $\Z$ with respect to the shift action of $\Z$ on $\Zmod{2}^{\oplus \Z}$. In particular, torsion subgroups of the lamplighter group have orders which are powers of $2$.

Shortly afterwards W. Dicks and T. Schick described in \cite{Dicks_Schick} an operator $T$ from the group ring of $(\Zmod{2}\wr \Z)^2$ and an heuristic evidence on why $\dimvn \ker T$ is probably irrational. Nonetheless, the question of irrationality of that specific number has remained open.

A breakthrough came in 2009, when T. Austin showed the following theorem.

\begin{theorem*}[\cite{arxiv:austin-2009}]
 The set of $l^2$-Betti numbers arising from finitely generated groups is uncountable.
\end{theorem*}
 
In particular there exist irrational $l^2$-Betti numbers. However, \cite{arxiv:austin-2009} did not provide a particular group which gives rise to irrational $l^2$-Betti numbers. 

\subsec{Results}

\begin{theorem}\label{thm_lamplighter}
The set of $l^2$-Betti numbers arising from the group $(\Zmod{2}\wr \Z)^3$ contains 
$$
\frac1{64}- \frac18\sum_{k=1}^\infty \frac1{2^{k^2+4k+6}},
$$
which is an irrational number.
\end{theorem}

Irrationality of the number above follows from the fact that its binary expansion is non-periodic. The author does not know whether the current transcendence results cover this number.

It is of some theoretical interest to have explicit finitely presented examples, so we point out also the following corollary. 

\begin{cory}\label{cory_lamplighter}
 Let $G$ be a group given by the presentation 
$$
	\langle a,t,s\,|\,a^2=1,[t,s]=1,[t^{-1}at,a]=1,s^{-1}as=at^{-1}at\rangle.
$$ 
The set of $l^2$-Betti numbers arising from $G^3$ contains 
$$
\frac1{64}- \frac18\sum_{k=1}^\infty \frac1{2^{k^2+4k+6}}.
$$
\end{cory}

In both Theorem \ref{thm_lamplighter} and Corollary \ref{cory_lamplighter}, the appropriate matrix whose kernel dimension is as stated can be explicitly written down in terms of the standard generators.

\begin{theorem}\label{thm_finitely_generated}
The set of $l^2$-Betti numbers arising from finitely generated groups is equal to the set of non-negative real numbers.
\end{theorem}

The group which realizes a given real number $r$ is ``as explicit as the binary expansion of $r$``. 

In other words, measures of atoms in spectral measures of rational group ring elements are arbitrary real numbers. The atoms themselves however are, at least for sofic groups, algebraic numbers (see \cite{Thom:Sofic_groups_and_Diophantine_approximation}).

We can also say something about the set of $l^2$-Betti numbers arising from finitely presented groups. Recall that a set $\Si$ of natural numbers is called {\bf computable} if there exists a Turing machine which lists the elements of $\Si$ in the increasing order (in other words, there exists an algorithm which allows to compute subsequent elements of $\Si$).

We say that a real number $r$ has  {\bf computable binary expansion} if the fractional part of $r$ is of the form
$$
\sum_{i\in \Si} \frac{1}{2^i}
$$ 
for some computable set $\Si$.

\begin{theorem}\label{thm_finitely_presented}
The set of $l^2$-Betti numbers arising from finitely presented groups contains all numbers with computable binary expansions. 
\end{theorem}

Examples of numbers with computable binary expansion include all algebraic numbers, $\pi$ and $e$. A fairly well-known example of a number whose binary expansion is not computable is Chaitin's constant encoding the halting problem (see \cite{calude_chaitin:what_is_a_halting_probability}).

Theorem \ref{thm_finitely_generated} has been independently proven by M. Pichot, T. Schick and A. Żuk in \cite{arxiv:pichot_schick_zuk-2010}. They also proved a result similar to Theorem \ref{thm_finitely_presented}. Let us also mention two later developments: (1) In  \cite{arxiv:lehner_wagner-2010},  F. Lehner and S. Wagner show that $\call C(\Zmod{p} \wr F_d)$ contains irrational algebraic numbers, where $F_d$ is the free group on $d$ generators, $d>2$, $p\ge 2d-1$;  (2) In  \cite{arXiv:grabowski-2010-2} the present author shows that $\call C(\Zmod{p}\wr \Ga)$ contains transcendental numbers, for all $p>1$ and all groups $\Ga$ which contain an element of infinite order.  

\subsec{Summary}

Section \ref{section_groupoids} deals with discrete measured groupoids, which seem to be the right context to generalize known spectral computations of random walk operators on groups of the form $\wh X \rtimes \Ga$, where $\wh X$ is discrete abelian. In short, the spectrum of an operator from a groupoid ring is shown to be the expected spectral measure of the convolution operator on a random (groupoid) Schreier diagram associated to the operator. The main computational tools are Proposition \ref{prop_computational_tool} and Corollary \ref{cory_computational_tool}. Section \ref{section_groupoids_examples} explains how they generalize a computation of R. Grigorchuk and A. Żuk from \cite{Grigorchuk_Zuk2001}, and a theorem of F. Lehner, M. Neuhauser and W. Woess from \cite{Lehner_Neuhauser_Woess:On_the_spectrum_of}.

Section \ref{section_tds1} introduces Turing dynamical systems, which are a natural generalization of Turing machines (possibly with oracles). The computational tool is used  in Theorem \ref{main_theorem} to show a connection between dynamical properties of a system and spectral properties of certain operator in the groupoid ring of the associated groupoid.

Section \ref{section_tds-examples} presents examples of Turing dynamical systems which are then used in Section \ref{section_Atiyah} to prove results about the Atiyah problem. 

\subsec{Open questions}

All questions are well-known to the experts.

\begin{question}
What is the set of $l^2$-Betti numbers arising from finitely presented groups?
\end{question}
Note that the set in question is countable. In \cite{arxiv:pichot_schick_zuk-2010} different numbers appear than those covered by Theorem \ref{thm_finitely_presented}. On the other hand, note that every $l^2$-Betti number is, by functional calculus, a limit of a  sequence $\tau_G( (1-T)^n)$, where $T \in M_k (\Q G)$ and $\tau_G$ is the group trace. If the group $G$ is sofic, then there are bounds known on what $n$ one has to take to be $\eps$-close to the limit (this follows from the determinant conjecture, see \cite{Elek_Szabo:Hyperlinearity_eseentially_free_actions_and}). It follows that if $G$ is a sofic group, then elements of $\call C(G)$ are computable by a Turing machine with an oracle for the word problem of $G$. If $G$ is finitely presented then the word problem of $G$ is known to be not harder than the halting problem. This gives some bound on what $l^2$-Betti numbers can arise from (sofic) finitely presented groups; however this bound seems to be far away for the techniques presented here 
or in \cite{arxiv:pichot_schick_zuk-2010}.

\begin{question}
 For a group $G$ and a ring $k\subset \C$, define $\call C(G, k)$ to be the set of those $r$ such that there exists $T\in kG$ with $\dimvn \ker T = r$. By definition, $\call C(G, \Q) = \call C(G)$. Is it true that for every group $G$ we have $\call C(G) = \call C(G, \C)$? In particular, is it true that the set $\call C(G, \C)$ is countable?
\end{question}
The answers are trivially yes for those torsion-free groups (or groups with bounded torsion subgroups) for which the so called \textit{strong Atiyah conjecture} holds: $\call C(G, \C) = \N$ (this has to be modified appropriately for bounded torsion groups). Examples include free groups and bounded-torsion elementary amenable groups (see \cite{Linnel:Atiyah_conjecture_for_free_groups}). On the other hand, the answers are not known even for $\Zmod{2}\wr\Z$. This motivates the next question.

\begin{question}
 What is $\call C(\Zmod{2}\wr \Z)$?
\end{question}
In Theorem \ref{thm_lamplighter} we prove $\call C ((\Zmod{2}\wr \Z)^3)\nsubseteq \Q$. After the first version of this article was submitted to arXiv, F. Lehner and S. Wagner showed in \cite{arxiv:lehner_wagner-2010} that $\call C(\Zmod{p} \wr F_d)$ contains irrational algebraic numbers, where $F_d$ is the free group on $d$ generators, $d>2$, $p\ge 2d-1$, which subsequently led the author to show in \cite{arXiv:grabowski-2010-2} that $\call C(\Zmod{p}\wr \Z)$ contains transcendental numbers, for all $p>1$. This raises the question whether $\call C(\Zmod{p}\wr \Z)$ contains irrational algebraic numbers. In fact, $\call C((\Zmod{p}\wr \Z)^k)$ contains algebraic numbers of degree $k$, as will be shown in a future version of \cite{arXiv:grabowski-2010-2}.

\begin{question}
Are $l^2$-Betti numbers \textbf{of} a countable discrete group rational?
\end{question}
For a precise definition of $l^2$-Betti numbers of a group, see e.g. \cite{Eckmann_intro} or \cite{Lueck:Big_book} for a more general definition. If $r$ is an $l^2$-Betti number of a group $G$, then it is in particular an $l^2$-Betti number arising from $G$, but not the other way around. All the examples in the literature so far, of groups which give rise to irrational $l^2$-Betti number, have an infinite normal amenable subgroup. This implies that all their $l^2$-Betti numbers are $0$ (see \cite{Lueck:Big_book}, Theorem 7.2).

\subsec{Thanks and acknowledgments}

The article is a modified version of author's thesis, defended on 10.03.2011 at Georg-August-Universität in Göttingen, and written under the advice and supervision of Andreas Thom and Thomas Schick. Parts of the manuscript had been revised while the author had been supported by EPSRC at Imperial College London and Oxford University.

The article is a continuation of the ideas contained in the papers \cite{Grigorchuk_Zuk2001},
\cite{Dicks_Schick}, and \cite{arxiv:austin-2009}. In particular, the main innovation  - Turing
dynamical systems - was largely motivated by the ``pattern recognition'' idea of
\cite{arxiv:austin-2009}. On the other hand, Section \ref{section_groupoids} came about as an attempt to understand the proof of Lemma 3.6 from \cite{Dicks_Schick}.

I thank Thomas Schick and Andreas Thom for helpful discussions; in particular for explaining to
me the aforementioned ``pattern recognition'' idea of Tim Austin. 

I owe special thanks to Andreas Thom, for encouragements to write down early versions of this work, and to Manuel Koehler, for committing many hours to discussing details of Section \ref{section_groupoids}.

I have had a great luck to have very critical readers of the draft versions of this article. Comments of Robin Deeley, 
Światosław Gal, Jarek Kędra, Thomas Schick and  Andreas Thom greatly improved the exposition, as
well as eliminated a number of typos and grammatical mistakes.

I also thank Franz Lehner, Stefaan Vaes and anonymous referee for helpful comments.

Finally,  parts of the work presented here were completed during my stay at Université
de la Méditerranée in Marseille in November 2009. My visit there was arranged thanks to  Michael Puschnigg, whom I
would like to whole-heartedly thank for this and for his kind hospitality.

\subsec{Notation, conventions and a lemma}\label{subsec_notation}

The term {\bf measurable space} refers to a standard Borel space. The word subset means measurable subset, whenever it makes sense; similarly for functions. We do not always explicitely check that the sets and functions we work with are measurable, but in all the cases such checks are standard. Axiom of Choice is not used. 

The cardinality of a set $S$ is denoted by $|S|$. The Hilbert space whose orthonormal basis consists of elements of $S$ is denoted by $l^2(S)$ or $l^2S$. The standard basis vectors are denoted by $\zeta_s$, $s\in S$. If $U$ is a subset of $S$ then $\chi_U\colon S \to \{0,1\}$ and $\chi(U)\colon S \to \{0,1\}$ denote the characteristic function of $U$. 

A sequence $(\ldots, x_{i-1}, x_i, x_{i+1}, \ldots)$ is denoted by $(x_i)$ or $\dot x_i$.

\paragraph{Graphs and diagrams} If $g$ is a graph then its sets of vertices and edges are denoted by $V(g)$ and $E(g)$. If $g$ is an oriented graph and $e\in
E(g)$, then the starting and ending points of $e$ are denoted by $s(e)$ and $r(e)$ ($r$ stands for \textit{range}).
The Hilbert space $l^2(V(g))$ is denoted also by $l^2g$. 

A full subgraph of a graph $g$ is a subgraph $h$ with the property that if $e\in E(g)$ is an edge between two vertices in $V(h)$, then $e\in E(h)$. 

A vertex in an oriented graph with no incoming edges, perhaps with an exception of a self-loop, is called a {\bf starting vertex}. Similarly for  a {\bf final vertex}.

When $A$ is a set or a sequence, an {\bf $A$-diagram} is an oriented graph of bounded degree with edge labels in $A$. $\C$-diagrams are called simply {\bf diagrams}. Given a diagram $g$, the associated {\bf convolution operator} is the unique operator $T\colon l^2g \to l^2g$ such that $\langle T(\zeta_x), \zeta_y\rangle$ is equal to $0$ if there are no edges between $x$ and $y$ and to sum of all the labels between $x$ and $y$ otherwise. If $g$ is a rooted diagram then the rooted spectral measure of the convolution operator $T_g$ on $g$ is the usual projection-valued spectral measure composed with the functional $P \mapsto \langle P \zeta_x, \zeta_x\rangle$, where $x$ is the root of $g$. 

If $A$ is a sequence, then the labels are assumed to ''remember`` the order of elements of $A$, e.g. if $A$ is the sequence $(a_1, a_2)$ then the two $A$-diagrams consisting of one oriented edge with label $a_1$ or $a_2$ are not the same $A$-diagram, even if $a_1=a_2$.

Given an oriented graph $g$, we use the expression ``diagram $g$`` for the diagram which has $g$ as the underlying graph and all the labels are $1$. If $g$ is an $A$-diagram we say ''graph $g$`` for the underlying oriented graph of $g$.

\paragraph{Rings and operators} Given a ring $R$, the ring of $k\times k$-matrices over $R$ is denoted by $M_k(R)$. A {\bf trace} $\tau$ on $R$ is a function $\tau\colon R \to \C$ such that $\tau(ab)=\tau(ba)$. The standard trace (i.e. sum of diagonal elements) on $M_k(\C)$ is denoted by $\tr$. If $R$ is a ${}^*$-ring of operators on a Hilbert space then we also require that $\tau(T^*T)$ is a non-negative real number, for all $T\in R$. If $R$ is an algebra over a field $\bb F \subset \C$ and $\tau$ is a trace on $R$, then we have an {\bf induced trace} on $M_k(R) \cong M_k(\bb F) \otimes R$ given by $\tr\otimes \tau$, The induced trace is also denoted by $\tau$.

If $R$ is a ${}^*$-ring of operators on a Hilbert space, then a trace $\tau$ on $R$ is called {\bf normal} if it extends to a continuous trace on the weak completion  of $R$. A trace is {\bf faithful} if, for every $T$, $\tau(T^*T)=0$ implies $T=0$. All traces we will consider are faithful and normal.

If $R$ is a ${}^*$-ring of operators  on a Hilbert space, $\tau$ is a faithful normal trace on $R$, and $T^*=T\in R$ then the {\bf spectral measure of $T$} is the usual projection-valued spectral measure of $T$ composed with $\tau$ (it makes sense to evaluate $\tau$ on spectral projections of $T$, since the latter are in the weak completion of $R$). The spectral measure of the set $\{0\}$ is called  {\bf von Neumann dimension of kernel of $T$}, denoted by $\dimvn \ker (T)$. For a non-self-adjoint $T$ we define $\dimvn \ker (T) = \dimvn \ker (T^*T)$.

We say that the spectral measure of $T^*=T \in R$ is {\bf pure-point}, or that $T$ has pure-point spectrum, if the spectral measure of $T$ is a countable sum of measures supported on single points.

\begin{lemma}\label{lemma_trace_preserving_homo}
 Let $A$ and $B$ be ${}^*$-rings of operators Hilbert spaces $H_A$ and $H_B$ with faithful normal traces $\tau_A$ and $\tau_B$. Let $\phi: A\to B$ be a trace preserving ${}^*$-homomorphism and $T^*=T\in A$. Then the spectral measure of $T$ with respect to $\tau_A$ is the same as the spectral measure of $\phi(T)$ with respect to $\tau_B$.
\end{lemma}

\begin{proof}
Since $\phi$ is ${}^*$-preserving, $\phi(T)$ is also self-adjoint. The spectral measure, as any $\si$-additive measure, is determined by measures of intervals. Let $I$ be an interval and $p_n$ be a sequence of polynomials converging to $\chi_I$ pointwise, everywhere on $\R$. By the definition of the spectral measure, we need to show $\tau_A(\chi_I(T)) = \tau_B(\chi_I(\phi(T))$.

By the spectral theorem, we have $p_n(T) \to \chi_I(T)$, and $p_n(\phi(T)) \to \chi_I(\phi(T))$.  Since $\phi$ is a homomorphism, we have $p_n(\phi(T)) = \phi(p_n(T))$. Since $\tau$ is normal, we have $\tau_A(p_n(T)) \to \tau_A(\chi_I(T))$, and $\tau_B(\phi(p_n(T))) \to \tau_B (\chi_I(\phi(T))$. The claim follows since $\phi$ is assumed to be trace preserving, in particular $\tau_B(\phi(p_n(T))) = \tau_A(p_n(T))$.
\end{proof}

\section{Groupoids}\label{section_groupoids}

\subsec{Definitions}

For more information on groupoids see \cite{Sauer_Thom:A_Spectral_sequence_to_compute} and references therein. 

A {\bf groupoid} $\call G$ is a small category whose morphisms are all invertible. The sets of objects and morphisms are denoted respectively by $\call G_0$ and $\call G$. The embedding $\bf 1\colon\call G_0 \to \call G$ sends an object $x$ to the identity morphism on $x$. The space $\call G_0$ will be often identified with a subset of $\call G$ via this embedding.

The maps $s,r\colon\call G \to \call G_0$, source and range maps, associate to a morphism its domain and codomain. Composition is a map $\call G {}_r\!\times_s \call G \to \call G$; composition of morphisms $\ga\colon x\to y$ and $\ga'\colon y\to z$ is denoted by $\ga\ga'$ ot $\ga\bullet \ga'$. Given $\ga\colon x\to y$, the inverse of $\ga$ is denoted either by $i(\ga)$ or by $\ga^{-1}$. Note that the symbol $\circ$ stays reserved for the standard composition of functions. 

For  $x\in \call G_0$, the sets $s^{-1}(x)$ and $r^{-1}(x)$ are denoted by $s^*x$ and $r^*x$. The set of those objects $y$ for which there exists a morphism between $x$ and $y$ is  the {\bf orbit of $x$}, denoted by $\call Gx$.

A {\bf discrete measurable groupoid} is a groupoid together with a structure of a measurable space on $\call G$, and such that $\call G_0$ is a measurable subset, fibers of the maps $s$ and $r$ are countable, and  the structure maps $s$, $r$, $i$ and composition are measurable. 

A {\bf discrete measured groupoid} is a discrete measurable groupoid $\call G$ together with a measure $\mu$ on $\call G_0$, such that the measures 
$$
\call G \supset U \mapsto \int_{\call G_0} |r^{-1}(x)\cap U| d\mu(x)
$$
and 
$$
\call G \supset U \mapsto \int_{\call G_0} |s^{-1}(x)\cap U| d\mu(x)
$$
are equal. This measure on $\call G$ is also denoted by $\mu$. 

From now on all groupoids will be discrete measured, unless explicitly stated otherwise. The following lemma is a direct consequence  of the definition of measure on $\call G$.

\begin{lemma}\label{lemma_injective_projection}
Suppose $U\subset \call G$ is such that $r$ restricted to $U$ is an injection. Then the measure of $U$ in $\call G$ is the same as the measure of $r(U)$ in $\call G_0$. 
\end{lemma}

For $x \in \call G_0$ let {\bf stabilizer $\st(x,\call G)$ of $x$} be the group of those $\al\in \call G$ such that $s(\al)=r(\al)=x$ We say that a groupoid $\call G$ is a {\bf relation groupoid}, if for almost all $x$ we have  $|\st(x, \call G)|=1$. If $\call G$ is a relation groupoid then we freely use the identification of $s^*x$ with $\call Gx$ given by $s^*x \ni \ga \mapsto r(\ga)$.

\begin{lemma}\label{lemma_integration_change_A_and_B}
Let $\call G$ be a groupoid, $A, B \subset \call G_0$, and $f\colon \call G_0 \to \C$ be constant on orbits of $\call G$. Then
$$
  \int_A f(x)\cdot |\st(x,\call G)|\cdot |\call Gx \cap B| \,d\mu(x) = \int_B f(x)\cdot |\st(x,\call G)|\cdot |\call Gx \cap A| \,d\mu(x)
$$
\end{lemma}
\begin{proof}
Let $U\subset \call G$ denote the set $\{\ga \in \call G\colon s(\ga)\in A, r(\ga)\in B\}\subset \call G$, and $\wt f:= f\circ s$. Both quantities are equal to
$$
  \int_{U} \wt f(x) \,d\mu(x)
$$
\end{proof}

A {\bf measurable edge} is a pair $(U, \phi)$, where $U\subset \call G_0$ and $\phi: U \to \call G$, such that $s\circ \phi\colon  U\to \call G_0$ is the identity embedding, and $r\circ \phi\colon  U \to \call G_0$ is injective. Note that $\phi$ and $r\circ \phi$ are automatically measure preserving. For the most part, we  write simply $\phi$, with the understanding that $U = \Dom( \phi)$ is the domain of definition of $\phi$. If $\phi$ is a measurable edge, then $\phi^{-1}$, the {\bf  inverse of $\phi$}, is the measurable edge with $\Dom (\phi^{-1}) = r\big(\Im(\phi)\big)$, and such that $\phi^{-1}\Big( r\big(\phi(x)\big) \Big) = i(\phi(x))$.

\begin{lemma}\label{lemma_countable_family_of_edges}
 There exists a countable family of measurable edges whose images are disjoint and such that the union of all their images is all of $\call G$.
\end{lemma}
\begin{proof}
The statement follows from a theorem of Luzin and Novikov (see \cite{Kechris:classical_descriptive_set_theory}, Theorem 18.10): there exists a division of $\call G$ into countably many disjoint sets such that the restriction of $s$ to any of them is injective (this is true for any measurable map with countable fibers). 
\end{proof}

Let $\Phi = (\phi_1, \ldots, \phi_n)$ be a sequence of measurable edges. The {\bf Schreier diagram} $\schr_\Phi(\call G)$ of $\call G$ with respect to $\Phi$ is the $\Phi$-diagram whose vertices are points of $\call G_0$, and there is an oriented edge between $x$ and $y$ with label $\phi_i$ iff there is a morphism between $x$ and $y$ in the image of $\phi_i$. Similarly the {\bf Cayley diagram} $\cay_\Phi(\call G)$ is the $\Phi$-diagram whose vertices are points of $\call G$, and there is an oriented edge between $\al$ and $\be$ with label $\phi_i$ if there is a morphism $\ga$ in the image of $\phi_i$ such that $\al\ga=\be$.

The full subdiagrams of $\schr_\Phi(\call G)$ and $\cay_\Phi(\call G)$ whose vertices are respectively $\call Gx$, and $s^*x$, where $x\in \call G_0$, are denoted by $\schr_\Phi(x, \call G)$ and $\cay_\Phi(x, \call G)$, and called {\bf Schreier and Cayley diagrams of $x$}. Both $\schr_\Phi(x, \call G)$ and $\cay_\Phi(x, \call G)$ are often considered as being rooted with root being the point $x$. If $\call G$ is a relation groupoid then $\schr_\Phi(x, \call G)$ and $\cay_\Phi(x, \call G)$ are naturally identified with each other.  The underlying oriented graphs of $\schr_\Phi(x, \call G)$ and $\cay_\Phi(x, \call G)$ are called {\bf Schreier and Cayley graphs} and they are denoted by the same symbols.

Let $T$ be a finite sum $\sum_{i=1}^n \phi_i f_i$ (for now purely formal), where $\phi_i$ are measurable edges and $f_i\in L^\infty(\call G_0)$, and let $\Phi$ be the sequence $(\phi_1, \ldots, \phi_n)$. The $\C$-diagram $\cay_T(x, \call G)$ is obtained from the $\Phi$-diagram $\cay_\Phi(x, \call G)$ by replacing each edge labeled $\phi_i$ starting at a vertex $\al\in \call G$ by an edge labeled $f_i\circ r(\al)$. Similarly for $\schr_T(x, \call G)$ we replace each edge labeled $\phi_i$ starting at a vertex $x\in \call G_0$ by an edge labeled $f_i(x)$.


\subsec{Groupoid ring}

Given a measurable edge $\phi$, define a bounded operator on $L^2 \call G $, denoted also by $\phi$, by
\begin{equation}\label{eq_def_action}
 \phi (F)(\ga) =  \left\{  
	\begin{array}{l l}
  		F\big(\ga\bullet (\phi^{-1}\circ r)(\ga)\big)  &\quad \mbox{if $r(\ga)\in \Dom (\phi^{-1})$,}\\
		0 & \quad \mbox{otherwise,}\\ 
	\end{array} \right. 
\end{equation}
where $F\in L^2(\call G)$.

Given $f\in L^\infty (\call G_0)$ and $F\in L^2(\call G)$, we define $f(F)\in L^2(\call G)$ to be $f(F)(\al):= (f\circ r)(\al) \cdot F(\al)$. This is an action of $L^\infty(\call G_0)$ on $L^2(\call G)$. 

The groupoid ring of $\call G$, denoted by $\C \call G$, is the ring of bounded operators on $L^2(\call G)$ generated by measurable edges and $L^\infty(\call G_0)$. By Lemma \ref{lemma_groupoid_ring_basics} below, $\C\call G$ is ${}^*$-closed.

Given a measurable edge $\phi$ and a set $U\subset \Dom (\phi)$, we let $\phi_{|U}$ and $\phi_U$ denote the restriction of $\phi$ to $U$. 

\begin{lemma}\label{lemma_groupoid_ring_basics}
 Let $\phi$ be a measurable edge and $f\in L^\infty (\call G_0)$. Then, in $\C G$, we have
\begin{enumerate}
 \item $\phi f = \phi_{|\supp (f)\cap \Dom (\phi)} f$,
 \item If $\Dom(\phi) \subset U$, and $\chi$ is the characteristic function of $U$ then $\phi \chi = \phi$.
 \item $\phi^* = \phi^{-1}$,\label{lemma_groupoid_ring_basics_star}
 \item $\phi f = \phi(f) \phi,$\label{lemma_groupoid_ring_basics_presentation}
\end{enumerate}
where $\phi(f)\in L^\infty(\call G_0)$ is defined by the formula
\begin{equation}\label{eq_def_action2}
\phi(f)(x) = \left\{ 
	\begin{array}{l l}
  		 (f\circ r\circ\phi^{-1})(x)  &\quad \mbox{if $x\in \Dom (\phi^{-1})$,}\\
		0 & \quad \mbox{otherwise.}\\ 
	\end{array} \right. 
\end{equation}
\end{lemma}
\begin{proof}
We only prove (4). The other statements are proved similarly and left to the reader as an exercise.

Let $F\in L^2 (\call G)$. Then $f(F) (\al) = f\circ r( \al)\cdot  F (\al)$ and therefore $\big[(\phi f)(F)\big] (\al)$ equals
$$
\begin{array}{l l}
  		 (f\circ r)\big(\al \bullet (\phi^{-1}\circ r)(\al)\big)\cdot F \big(\al\bullet(\phi^{-1}\circ r)(\al)\big)  & \quad\mbox{if $r(\al) \in \Dom(\phi^{-1})$,}\\
		0 &  \quad\mbox{otherwise.}\\
\end{array}
$$
From \eqref{eq_def_action} we get
$$
\big[(\phi (f)\phi)(F)\big] (\al) =\left\{ 
	\begin{array}{l l}
  		 \big(\phi(f)\circ r\big)(\al) \cdot F (\al \bullet (\phi^{-1}\circ r)(\al))  &\quad \mbox{if $r(\al) \in \Dom(\phi^{-1})$,}\\
		0 & \quad \mbox{otherwise,}\\ 
	\end{array} \right. 
$$
but, from \eqref{eq_def_action2} we see that, if $r(\al) \in \Dom(\phi^{-1})$ then  $(\phi(f)\circ r)(\al) = (f\circ r\circ \phi^{-1}\circ r)(\al)$, and so the claim follows from noting that $r(\al\be)=r(\be)$ for every composable pair $\al, \be$ of morphisms. 
\end{proof}

In particular, each element of $\C \call G$ can be (non-uniquely) represented by a finite linear combination of operators $\phi\cdot f$, where $f\in L^\infty (\call G_0)$, and $\phi$ is a measurable edge. If $T\in \call G$ is represented by $\sum \phi_i f_i$ then we abuse notation by denoting the sum $\sum \phi_i f_i$ also by $T$.

The trace $\tau_{\call G}$ on $\C \call G$ is defined by the formula
$$
 \tau_ {\call G}(T) := \langle T\chi_0, \chi_0 \rangle_{L^2 \call G}, 
$$
where $\chi_0$ is the characteristic function of $\call G_0 \subset \call G$. It is positive, faithful and normal.

\subsec{Action groupoids and Pontryagin duality}

Let $\Ga$ be a discrete countable group, $(X,\mu)$ be a probability measure space and $\rho\colon \Ga \actson  X$ be a right measure preserving action, which is not necessarily free. The {\bf action groupoid $\call G(\rho)$} is a measured groupoid whose space of objects is $X$, and whose space of morphisms is $X\times \Ga$. The structure maps are given by $s(x,\ga) =x$, $r(x,\ga) = \rho(\ga)(x)$. Composition of $(x, \al)$ and $(\rho(\al)(x), \be)$ is $(x, \al\be)$. The inverse of $(x,\ga)$ is $(\rho(\ga)(x), \ga^{-1})$.

Each element $\ga \in \Ga$ gives rise to a measurable edge $x\mapsto (x,\ga)$, denoted also by $\ga$, whose domain of definition is all of $X$.

For the rest of this subsection, $(X, \mu)$ is a compact abelian group with the normalized Haar measure and the action $\rho\colon \Ga \actson X$ is by continuous group automorphisms. The action $\wh \rho$ of $\Ga$  on the Pontryagin dual $\wh X$ of $X$ is defined as $\wh\rho (\ga) (f) (x) = f(\rho(\ga^{-1})(x))$, where $f\colon X \to \C$ is an element of $\wh X$. 

For more information on  Pontryagin duality see e.g. \cite{Folland:A_course_in_abstract_harmonic_analysis}. In particular, Pontryagin duality induces a map $P\colon \wh X \to L^\infty (X)$. If $\wh x \in \wh X$ then $P(\wh x)$ is denoted also by $x$.

\begin{proposition}\label{prop_pontryagin}
There is a trace-preserving ${}^*$-embedding of the complex group ring $\C (\wh X \rtimes_{\wh \rho} \Ga)$ into the groupoid ring $\C \call G(\rho)$, which sends $ f\in \wh X$ to $P(f)$, and $\ga\in \Ga$ to a measurable edge $\ga$. This embedding will be denoted by $P\otimes 1$.

In particular, if $T=T^* \in \C (\wh X \rtimes_{\wh \rho} \Ga)$, then $T$ and $P\otimes 1(T)$ have the same spectral measures.
\end{proposition}
\begin{proof}
To begin with, we show that the map $\sum c_i \cdot \wh a_i \cdot \ga_i \mapsto \sum c_i\cdot a_i \cdot \ga_i$, $c_i \in \C$, $\wh a_i\in \wh X$, $\ga_i \in \Ga$, is a ring homomorphism. It is well-defined,  since every element of $\C (\wh X \rtimes_{\wh \rho} \Ga)$ can be written in a unique way as $\sum c_i \cdot \wh a_i \cdot \ga_i$. It certainly is a ring homomorphism when restricted to $\C\Ga$ and to $\C\wh X$. The standard presentation of a semi-direct product and Lemma \ref{lemma_groupoid_ring_basics}\eqref{lemma_groupoid_ring_basics_presentation} imply that it is a homomorphism on all of $\C (\wh X \rtimes_{\wh \rho} \Ga)$.

Lemma \ref{lemma_groupoid_ring_basics}\eqref{lemma_groupoid_ring_basics_star} implies that the ${}^*$-operation is preserved.

By linearity, it is enough to check the trace-preservation on an element of the form $\wh a \cdot \ga \in \C (\wh X \rtimes_{\wh \rho} \Ga)$, where $\wh a \in \wh X$, $\ga \in \Ga$. The group ring trace of $\wh a \cdot \ga$ is equal to $1$ if $\wh a$ and $\ga$ are the neutral elements of $\wh X$ and $\Ga$ respectively, and is equal to $0$ otherwise. 

We consider three cases.

(1) Both $\wh a$ and $\ga$ are the neutral elements. Then $a$ is the function on $X$ constantly equal to $1$ and $\tau_{\call G(\rho)}(P\otimes 1(\wh a \cdot \ga)) = \langle \chi_0, \chi_0 \rangle = 1$.

(2) $\ga$ is not the neutral element. Then the functions $\chi_0$ and $\ga(\chi_0)$ have disjoint supports and, therefore, also $\chi_0$ and $(a\ga)(\chi_0)$ have disjoint supports, so $\langle (a\ga)(\chi_0), \chi_0 \rangle = 0$.

(3) $\ga$ is the neutral element, but $\wh a$ is not. Then, it follows from Pontryagin duality that $a$ is a non-constant function on  $\call G_0$, and the trace we have to compute is equal to $\langle a, \chi_0 \rangle = \int a\, d\mu$. Let $x, y\in \call G_0$ with $a(x) \neq a(y)$. We use that $a$ is a group homomorphism and invariance of Haar measure to get
$$
\int a(z) \, d\mu(z) = \int a(xz) \, d\mu(z) = a(x)\int a(z) \, d\mu(z).
$$
Repeating with $y$ in place of $x$ we obtain $a(x)\int a(z) \, d\mu(z)=a(y)\int a(z) \, d\mu(z)$, which is possible only if $\int a(z) \, d\mu(z) = 0$.

The statement about the spectral measures follows from Lemma \ref{lemma_trace_preserving_homo}.
\end{proof}

\begin{lemma}\label{lemma_cylinder_sets}
Let $\Z(\frac{1}{2})$ be the subring of $\Q$ generated by $\Z$ and $\frac{1}{2}$. If $X$ is a product of infinitely many copies of $\Zmod{2}$ indexed by a set $I$, then the  image of $\Z(\frac{1}{2}) (\wh X \rtimes_{\wh \rho} \Ga)$ under $P\otimes 1$ is generated over $\Z(\frac{1}{2})$ by measurable edges $\ga\in \Ga$ and the characteristic functions of cylinder sets.
\end{lemma}
\begin{proof}
 Let $R\subset \C \call G$ be the ring generated by characteristic functions of cylinder sets and measurable edges $\ga\in \Ga$.

First, we show that image of $\Z(\frac{1}{2}) (\wh X \rtimes_{\wh \rho} \Ga)$ is contained in $R$. Clearly, $(P\otimes 1)(\Ga) \subset R$. Note that $\wh X$ is a direct sum of infinitely many copies of $\Zmod{2}$ indexed by $I$. Let $g_i$ be the generator of $\Zmod{2}$ corresponding to the index $i\in I$. Direct computation shows that $P\otimes 1( \frac{e+g_i}{2})$ is the characteristic function of the cylinder set $\{(x_j)\in \prod_I \Zmod{2}\colon x_i = 0\}$. Also $P\otimes 1 (e)$ is a characteristic function of a cylinder set (namely, of the whole $X$). The statement follows, since $\Z(\frac{1}{2}) (\wh X)$ is generated, as a $\Z(\frac{1}{2})$-ring, by $\frac{e+g_i}{2}$ and $e$.

In the other direction, we just saw that the characteristic functions of cylinder sets $\{(x_j)\in \prod_I \Zmod{2}\colon x_i = 0\}$ are in the image. Since the constant function $1$ is also in the image, it follows that characteristic function  of $\{(x_j)\in \prod_I \Zmod{2}\colon x_i = 1\}$ is in the image as well. Every cylinder set is an intersection of sets of those two types, so the claim follows.
\end{proof}

\subsec{Groupoid ring and convolution operators on random diagrams}

For definitions and notation on direct integrals of spaces and operators see \cite{Folland:A_course_in_abstract_harmonic_analysis}, Chapter 7.4. Unless explicitly stated otherwise, all integrals are taken over the space $\call G_0$.

Consider the field $x \mapsto l^2(s^* x)$ of Hilbert spaces over $\call G_0$.  For a measurable edge $\phi$, define a section $S_\phi$, by
$$
  S_\phi(x) = \left\{ 
	\begin{array}{l l}
  		0 & \quad \mbox{if $x\notin \Dom (\phi)$}\\
		\zeta_{\phi(x)} & \quad \mbox{otherwise}\\ 
	\end{array} \right. 
$$
 Let $\psi_i$ be a countable family of measurable edges from Lemma \ref{lemma_countable_family_of_edges}. We make $l^2(s^*x)$ into a measurable field by equipping it with the family of sections $S_{\psi_i}$.

Given $T\in \C \call G$ represented by a finite sum $\sum_{i=1}^n \phi_i \cdot f_i$, we define a field of operators $T_x  \colon l^2 (s^*x) \to l^2 (s^*x)$ by setting $T_x$ to be the convolution operator on the diagram $\cay_T(x, \call G)$. The following proposition says that  the spectral measure of $T$ is equal to the expected rooted spectral measure of the convolution operator $T_x$ on the random diagram $\cay_T(x,\call G)$.

\begin{proposition}\label{prop_direct_integral_isomorphism}
There is an isomorphism of Hilbert spaces $L^2(\call G)$ and 
$$
\int^\oplus l^2( s^*x) \, d\mu(x)
$$ 
which sends a function $F\colon \call G \to \C$ to a section
\begin{equation}\label{eq_direct_integral_isomorphism}
x \mapsto \sum_{\ga \in s^*x} F(\ga)\cdot \zeta_\ga.
\end{equation}
Under this isomorphism elements of $\C \call G$ are decomposable. Operator $T\in \C \call G$ corresponds to the operator $\int^\oplus T_x \, d\mu(x)$. Furthermore, 
$$
  \tau_{\call G}(T) = \int \langle T_x( \zeta_x), \zeta_x\rangle_{l^2(s^*x)} \, d\mu(x).
$$
\end{proposition}
\begin{proof}
The statements about the decomposition of $T$ and the trace follow from the formula \eqref{eq_direct_integral_isomorphism} through a direct computation.

We need to check \textit{(i)} that the formula \eqref{eq_direct_integral_isomorphism} defines a measurable element of the field  $\int^\oplus l^2( s^*x) \, d\mu(x)$, \textit{(ii)} that this field is square-summable, and that the resulting map of Hilbert spaces is \textit{(iii)} isometric and \textit{(iv)} surjective.

\textit{(i)} For each measurable edge $S_\phi$ we need to check that the function 
$$
x\mapsto \langle S_\phi(x), \sum_{\ga \in s^*x} F(\ga)\cdot \zeta_\ga \rangle_{l^2s^*x}
$$ 
is measurable. By definition of $S_\phi$, this function is non-zero only on $\Dom(\phi)$, where it is equal to
$$
\langle \zeta_{\phi(x)}, \sum_{\ga \in s^*x} F(\ga)\cdot \zeta_\ga \rangle_{l^2s^*x},
$$
which is equal to $F(\phi(x))$.

\textit{(ii) and (iii)} Due to Lemma \ref{lemma_countable_family_of_edges}, we have 
$$
\langle F, F \rangle_{L^2 \call G} = \int_{\call G}  |F(\ga)|^2 \, d\mu(\ga) = \int  \sum_{\ga \in s^*x} | F(\ga)|^2 \,d\mu(x),
$$
which is equal to
$$
\int \langle \sum_{\ga \in s^*x} F(\ga)\cdot \zeta_\ga, \sum_{\ga \in s^*x} F(\ga)\cdot \zeta_\ga \rangle_{l^2s^*x} \, d\mu(x).
$$

\textit{(iv)} Let an element of $\int^\oplus l^2( s^*x) \, d\mu(x)$ be given by a measurable section $F(x)\in l^2(s^*x)$. Define $F \in L^2(\call G_0)$ as $F(\ga) = \langle F(s(\ga)), \zeta_\ga \rangle$. By \eqref{eq_direct_integral_isomorphism}, the image of $F$ is the section $F(x)$.
\end{proof}

\subsec{Subgroupoids}

 Let $\Phi = (\phi_1, \ldots )$ be a sequence of measurable edges. The {\bf subgroupoid generated by $\Phi$}, denoted by $\call G (\Phi)$, is the discrete measured groupoid whose space of objects is equal to $\call G_0$ and whose morphisms are generated by all morphisms in the images of  $\phi_i$, $\phi_i^{-1}$, and by all identity morphisms.

A {\bf subgroupoid} of a groupoid is a subgroupoid generated by a sequence of measurable edges.

\begin{proposition}
 Let $\call H$ be a subgroupoid of $\call G$. Note that if $\phi$ is a measurable edge in $\call H$, then it is also a measurable edge in $\call G$. Similarly,  elements of $L^\infty  \call H_0$ are at the same time elements of $L ^\infty \call G_0$. These two identifications  extend to a trace-preserving ${}^*$-embedding $\C \call H \hookrightarrow \C \call G$.
\end{proposition}

\begin{proof}
Let $T\in \C \call H$ be represented by a finite sum $\sum_{i=0}^n \phi_i f_i$, and let $\Phi$ be the sequence $(\phi_1, \ldots, \phi_n)$. Let $S\in \C \call G$ be represented by the same finite sum. By Proposition \ref{prop_direct_integral_isomorphism} and the remark before, $\tau_{\call H} (T)$ is the expected rooted spectral measure of the convolution operator on the diagram $\cay_T(x, \call H(\Phi))$, and $\tau_{\call G} (S)$ is the expected rooted spectral measure of the convolution operator on the diagram $\cay_S(x, \call G(\Phi))$. But those two diagrams have isomorphic connected components of the root, and the rooted spectral measure depends only on the connected component of the root, and so $\tau_{\call H} (T)= \tau_{\call G} (S)$.

Therefore, we have a map from the set of finite sums representing elements of $\C \call H$ to the set of finite sums representing elements of $\C \call G$, which is trace preserving and ${}^*$-preserving. It follows that this map induces a well-defined linear ${}^*$-embedding of $\C \call H$ into $\C \call G$ because both $\tau_{\call H}$ and $\tau_{\call G}$ are faithful.
\end{proof}

\begin{cory} \label{cory_subgroupoids}
Let $T^*=T\in \C \call G$ be represented by a sum $\sum_{i=1}^n \phi_i f_i$, and let $\Phi=(\phi_1, \ldots, \phi_n)$. Then $T$ is in the image of the embedding $\C \call G (\Phi) \hookrightarrow \C \call G$, and the corresponding element of $\C \call \call G(\Phi)$ is also denoted by $T$. The spectral measure of $T$ in $\C \call \call G(\Phi)$ is the same as the spectral measure of $T$ in $\C \call G$.
\end{cory}

\begin{proof} Follows from Lemma \ref{lemma_trace_preserving_homo}
\end{proof}

\subsec{Finite groupoids}

A groupoid $\call G$ is {\bf finite} if for almost all points $x\in \call G_0$ the set $s^*x$ is finite. We say that a groupoid $\call G$ has {\bf finite orbits}, if almost all points in $\call G_0$ have finite orbits. Note that a relation groupoid with finite orbits is finite, and that a finite groupoid has finite orbits. 

In particular if $\call G$ is a finite groupoid, then there exists a {\bf fundamental domain}, i.e. a measurable subset $D\subset \call G_0$ such that every finite orbit intersects $D$ exactly once. 

\begin{proposition}\label{prop_computational_tool}
Let $\call G$ be a finite groupoid, and let $D$ be a fundamental domain of $\call G$. There is a  ${}^*$-representation of $\C \call G$ on $\int_D^\oplus l^2(s^*x) \, d\mu(x)$, which sends an operator $T \in \C \call G$ to 
$$
\int_D^\oplus T_x \, d\mu(x).
$$
Under this representation
$$
  \tau_{\call G} (T) = \int_D \frac{\tr(T_x)}{|\st(x,\call G)|} \, d\mu(x).
$$
In particular
$$
  \dimvn\ker T = \int_D \frac{\dim \ker T_x}{|\st(x, \call G)|} \,d\mu(x).
$$
\end{proposition}

\begin{proof}
 Let $D^c$ be the complement of $D$. We have the direct sum decomposition 
$$
\int^\oplus l^2(s^*x) \, d\mu(x) = \int_D^\oplus l^2(s^*x) \, d\mu(x) \oplus \int_{D^c}^\oplus l^2(s^*x) \, d\mu(x),
$$
and the corresponding decomposition of the operator $T$:
$$
\int T_x \, d\mu(x) = \int_D^\oplus T_x \, d\mu(x) \oplus \int_{D^c}^\oplus T_x \, d\mu(x).
$$
It follows that $T \mapsto \int_D^\oplus T_x$ is a ${}^*$-representation. Thus, we need only to show the statement about  the trace. The right hand side is equal to
$$
\int_D \sum_{\ga\in  s^*x} \frac{1}{|\st(x,\call G)|}\langle T_x\zeta_\ga, \zeta_\ga \rangle \, d\mu(x).
$$
Note that for all $x\in \call G_0$ and $\al\in G$ such that $r(\al)=x$ we have $\langle T_x \zeta_\ga, \zeta_\ga\rangle = \langle T_{s(x)}\zeta_{\al\ga}, \zeta_{\al\ga}$, therefore, by putting $\al = i(\ga)$, the above is equal to
$$
\int_D \sum_{\ga\in  s^*x} \frac{1}{|\st(x,\call G)|}\langle T_{r(\ga)}\zeta_{r(\ga)}, \zeta_{r(\ga)} \rangle \, d\mu(x),
$$
which is the same as 
$$
\int_D \sum_{y\in \call Gx} \langle T_{y}\zeta_y, \zeta_y \rangle \, d\mu(x),
$$
which, by Lemma \ref{lemma_injective_projection} and the definition of fundamental domain, equals
$$
\int_{\call G_0} \langle T_{x} \zeta_x, \zeta_x \rangle \, d\mu(x) = \tau_{\call G} (T).
$$
The ``in particular`` statement follows from Lemma \ref{lemma_trace_preserving_homo}.
\end{proof}

\begin{cory}\label{cory_computational_tool}\label{cory_computational_tool_neat}
Let $\call G$ be a finite groupoid and let $T\in \C\call G$. Then
$$
\tau_{\call G}(T) = \int \frac{1}{|\call Gx||\st(x,\call G)|} \tr(T_x) \,d\mu(x)= \int \frac{1}{|s^*x|} \tr(T_x) \,d\mu(x).
$$
\end{cory}

\begin{proof}
The second equality is clear. The first one follows from the proposition and Lemma \ref{lemma_integration_change_A_and_B} for $A=\call G_0$, $B=D$ and  $f(x)=\tr(T_x)$.
\end{proof}

\section{Examples}\label{section_groupoids_examples}

We show how Proposition \ref{prop_computational_tool} and Corollary \ref{cory_computational_tool} generalize known computations of spectra of random walk operators.

\subsec{Computation of Grigorchuk \& Żuk}\label{subsec_grigorchuk_zuk}

 We start by computing the von Neumann dimension of the kernel of a random walk operator on the group $\Zmod{2} \wr \Z$. This was originally done, by different methods, by R. Grigorchuk and A. Żuk in \cite{Grigorchuk_Zuk2001}. Compare also \cite{Dicks_Schick}. 

The lamplighter group $\Zmod{2} \wr \Z$ is defined as $\Zmod{2}^{\oplus \Z} \rtimes \Z$, where the action of $\Z$ on $\Zmod{2}^{\oplus \Z}$ is by the shift,  i.e. $[t((x_i))]_j = x_{j+1}$, 

Let $X=\Zmod{2}^\Z$. Given $\eps_{-k}, \ldots, \eps_l \in \Zmod{2}$, the set 
$$
\{(x_i)\in X\colon x_{-k} = \eps_{-k}, \ldots, x_l = \eps_l \}
$$ 
is denoted by $ [\eps_{-k}\,\ldots\,\underline{\eps_0}\,\ldots\,\eps_l]$ 
and its characteristic function is denoted by  
$$
\chi[\eps_{-k}\,\ldots\,\underline{\eps_0}\,\ldots\,\eps_l].
$$ 
Similarly the set $\{(x_i)\in X\colon x_{-1}=0\}$ is denoted by $[0\underline{\,\cdot\,}]$ and its characteristic function by 
$ \chi[0\underline{\,\cdot\,}]$. Concrete elements from the set $[\eps_{-k}\,\ldots\,\underline{\eps_0}\,\ldots\,\eps_l]$ are denoted by $(\eps_{-k}\,\ldots\,\underline{\eps_0}\,\ldots\,\eps_l)$.

\begin{theorem}[\cite{Grigorchuk_Zuk2001}] Let $T$ be the element in the rational group ring of the lamplighter group given by $T = \frac{1}{2}(t + t^{-1} + tg + gt^{-1})$, where $t$ is the generator of $\Z$, and $g \in \Zmod{2}^{\oplus \Z}$ is a characteristic function $\Z \to \Zmod{2}$  of $\{0\} \subset \Z$. Then $\dimvn\ker T = \frac{1}{3}$.
\end{theorem}

\begin{proof}
Let $X = \{0,1\}^\Z$ with the standard product measure, and let $\rho\colon \Z\, \actson X$ be the Bernoulli shift action, i.e. $[\rho(t)((x_i))]_j = x_{j+1}$, where $t$ is the distinguished generator of $\Z$. Let $\call G$ denote the action groupoid $\call G(\rho)$

Note that $T = t \cdot \frac{1+g}{2} + t^{-1}\cdot \frac{1+t(g)}{2}$, and that the Pontryagin dual of $\frac{1+g}{2}$ is $\chi[\,\underline 0\,]$ and of $\frac{1+t(g)}{2}$ is $\chi[0\underline{\,\cdot\,}]$. Therefore, by Proposition \ref{prop_pontryagin}, the spectral measure of $T$ is the same as the spectral measure of the operator  $t_{[\,\underline 0\,]} + t^{-1}_{[0\underline{\,\cdot\,}]}  \in \C \call G$. Let us call the latter operator $T$ as well and let $\Phi$ be the sequence $(t_{[\,\underline 0\,]}, t^{-1}_{[0\underline{\,\cdot\,}]})$ of measurable edges.

The set 
$$
  D := [1\underline 1] \cup \bigcup_{k>1} [1\underline 0 0^{k-1}1]
$$
is a fundamental domain for $\call G(\Phi)$. For $x\in (1\underline 0 0^{k-1}1)$ the Schreier graph $\schr_T(x, \call G(\Phi))$ is 
$$
\bullet\leftrightarrows \bullet \leftrightarrows \cdots \leftrightarrows \bullet \leftrightarrows \bullet
$$
with $k$ vertices. The kernel of the convolution operator  is $1$-dimensional if $k$ is odd and $0$-dimensional otherwise. For $x\in [1\underline]$, the Schreier graph consists of one vertex and no edges, and so the kernel of the convolution operators is $1$-dimensional. 

Applying Proposition \ref{prop_computational_tool} we get 
$$
\dimvn\ker T = \frac{1}{4}+\sum_{k>1,\ 2|k\,} \frac{1}{2^{k+2}} = \sum_{l=0}^\infty \frac{1}{2^{2l+2}} = \frac{1}{3}.
$$

\end{proof}

If we put more effort into computing the spectra of the convolution operators on the Schreier graphs in the proof, we could compute the whole spectral measure of $T$ (see \cite{Grigorchuk_Zuk2001} or \cite{Dicks_Schick}).

\subsec{Percolation theory, theorem of Lehner, Neuhauser \& Woess}\label{subsec_percolation}

We now present  a  theorem of F. Lehner, M. Neuhaser and W. Woess from \cite{Lehner_Neuhauser_Woess:On_the_spectrum_of}. 

Let $\Ga$ be a finitely generated group, $\Psi=(\ga^1,\ldots, \ga^n)$ be a  finite sequence of generators for $\Ga$ which is symmetric, i.e. for each $\ga\in \Ga$ the number of times $\ga$ appears in $\Psi$ is equal to the number of times $\ga^{-1}$ appears in $\Psi$. Let  $\call C$ be the associated directed Cayley graph. An {\bf animal} is a full connected subgraph of $\call C$ which contains the neutral element $e\in \Ga$. The {\bf root} of an animal is defined to be the neutral element $e\in \Ga$. For an animal $\call L$ let us put 
$$
\p_\Psi \call L := \{\al \in V(\call C)-V(\call L)\colon \exists \be\in V(\call L) \text{ such that there is an edge from $\al$ to $\be$}\}.
$$

Given a point $x \in \Zmod{p}^\Ga$ let $\call C_0(x)$ be the full subgraph of $\call C$ with $V(\call C_0(x))=\{\al\in \Ga\colon x(\al) =0\}$, and let $\call L(x)$ be the animal whose set of vertices is equal to the connected component of $\call C_0(x)$ containing $e\in \Ga$ (if $e\notin \call C_0G$ then $\call L(x)$ is the empty animal).

Let $T_{\call L}$ denote the (oriented) convolution operator $l^2 \call L \to l^2 \call L$.

For $q\in[0,1]$, let $\nu_q$ denote the product measure on $\Zmod{p}^\Ga$ of a fixed measure $\nu$ on $\Zmod{p}$ such that $\nu(\{0\}) = q$, $\nu(\{1,2,\ldots, p-1\})=1-q$. 
We say that a parameter $q$ is {\bf subcritical for $\Ga$ with respect to $\Psi$} if animals $\call L(x)$ are $\nu_q$-almost surely finite. Subcriticality does not depend on $p$ or a choice of $\nu$. For example, for the free group $F_k$ with the standard generating set every $q<\frac{1}{2k-1}$ is subcritical. See \cite{benjamini_schramm:percolation_beyond} for more information on percolation theory.

\begin{theorem}[\cite{Lehner_Neuhauser_Woess:On_the_spectrum_of}]
 Let $G = \Zmod{p} \wr \Ga$, let $\pi\in \Q\Zmod{p}$ be the projection $\pi := \frac1p\cdot \sum_{a\in \Zmod{p}} a$, let  $T\in \Q G$ be defined as
$$
T := \pi\cdot \left( \sum_{i=1}^n \ga^i \right) \cdot \pi.
$$
Suppose that the parameter $\frac{1}p$ is subcritical for $\Ga$ with respect to $\Psi$. Then the spectral measure of $T$ is pure-point and
$$
\dimvn \ker T = \sum_{\call L} \frac{1}{|\call L|}\cdot \left(\frac1p\right)^{|\call L|} \cdot \left(\frac{p-1}{p}\right)^{|\p_\Psi \call L|} \cdot \dim\ker T_{\call L},
$$
where the sum is over all finite animals.
\end{theorem}

\begin{proof}
Let $X = \Zmod{p}^\Ga$ with the normalized Haar measure $\mu$ and let the action $\rho\colon \Ga \actson X$ be the Bernoulli shift, i.e.
$$
[\rho(\ga)(x_\al)]_\be := x_{\be\ga}.
$$
The action groupoid $\call G(\rho)$ is denoted by $\call G$.

For $\ga\in \Ga$ let $[\underline 0\xrightarrow{\ga}0] \subset X$ be the set $\{f\in X\colon f(e) = 0, f(\ga)=0\}$ and $\chi[\underline 0\xrightarrow{\ga}0]$ be its characteristic function. Computation shows that the image of $T$ under the Pontryagin duality map $P\otimes 1\colon \Q (\Zmod{p} \wr \Ga ) \to \C \call G(\rho)$ is the element 
$$
\sum_{i=1}^n \ga^i \cdot \chi[\underline 0\xrightarrow{\ga}0],
$$
which we also denote by $T$. Let $\Phi$ be the sequence of $\ga_i$ of measurable edges defined as $\ga^i$ restricted to the set $[\underline 0\xrightarrow{\ga^i}0]$. For an animal $\call L$ define 
$$
X(\call L) = \{x\in X: \call L(x) = \call L\}.
$$
and note 
$$
  \mu (X(\call L)) = \left(\frac1p\right)^{|\call L|} \cdot \left(\frac{p-1}{p}\right)^{|\p_\Psi \call L|},
$$
since for a point $x$ being in $X(\call L)$ means $x(\al) = 0$ for $\al \in V(\call L)$ and $x(\al) \neq 0$ for $\al \in  \p_\Psi \call L$. The sets $X(\call L)$ are disjoint for different animals and by subcriticality
$$
 \mu (\bigcup_{\text{finite } \call L} X(\call L)) = 1.
$$
Computation shows that for a point $x\in X(\call L)$ the Schreier graph $\schr_\Phi(x,\call G(\Phi))$ is isomorphic to $\call L$. It follows that under this isomorphism $T_x$ corresponds to $T_{\call L}$, because both operators are the convolution operators. In particular the groupoid $\call G$ is finite, so Corollary \ref{cory_computational_tool} gives us 
$$
\dimvn\ker T = \int \frac{1}{|\call Gx|} \dim\ker T_x \,d\mu(x) = \sum_{\call L} \frac{1}{|\call L|}\cdot \int_{X(\call L)} \dim\ker T_{\call L} \,d\mu(x), 
$$
where the sum is over all finite animals. This shows the statement about $\dimvn\ker T$. 

The spectral measure of $T$ is pure-point because $T$ is similar to the operator $\oplus_{\call L} T_{\call L}$, were the sum is again over all finite animals.
\end{proof}

The author was informed by Franz Lehner that it is an open problem to determine whether $\frac{1}{p} > q_c$ implies that the continuous part appears in the spectral measure of the operator $T$ above, even for groups $\Ga=\Z^k$, $k>1$. Surprisingly to the author, a conjecture in mathematical physics states that for $k=2$ the continuous part does not appear.
\section{Turing dynamical systems}\label{section_tds1}

\subsec{Definitions and basic properties}

Let $(X, \mu)$ be a probability measure space and $\rho:\Ga\curvearrowright X$ be a right probability measure preserving action of a countable discrete group $\Ga$ on $X$.  A {\bf dynamical hardware} is the following data: $(X,\mu)$, the action $\rho$, and a division $X = \bigcup_{i=1}^n X_i$ of $X$ into disjoint measurable subsets. For brevity, we denote such a dynamical hardware by $(X)$.

Suppose now that we are given a dynamical hardware $(X)$ and we choose three additional distinguished disjoint subsets of $X$, each of which is a union of certain $X_i$'s: the initial set $I$, the rejecting set $R$, and the accepting set $A$ (all or some of them might be empty). Furthermore, suppose that for every set
$X_i$, we choose one element $\ga_i$ of the group $\Ga$ in such a way that the elements corresponding to the sets $X_i$ which are subsets of $R\cup A$ are equal to the neutral element  of $\Ga$. 

Define a map $T_X: X\to X$ by
$$
	T_X(x):= \rho(\ga_i)(x) \quad \text{ for } x\in X_i.
$$

A {\bf dynamical software} for a given dynamical hardware $(X)$ is the following data: the distinguished sets $I, A$ and $R$, the choice of elements $\ga_i$, and the map $T_X$, subject to the conditions above. The map $T_X$ will be called the {\bf Turing map}, and the whole dynamical software will be denoted by $(T_X)$. A {\bf Turing dynamical system} is a dynamical hardware $(X)$ together with a dynamical software $(T_X)$ for $(X)$.  We will denote such a Turing dynamical system by $(X, T_X)$.

\begin{proposition}\label{prop_measure_contracting}
 In any Turing system $(X,T_X)$, the Turing map $T_X$ is mea\-sure-con\-tra\-cting, i.e., for every measurable set $U\subset X$ the set $T_X(U)$ is measurable and $\mu(T_X(U)) \le \mu (U)$. If $T_X$ is injective on $U$ then $\mu(T_X(U)) = \mu(U)$.
\end{proposition}

\begin{proof}
$T_X(U)$ is measurable since $T_X$ is finite-to-one and measurable.

Define $U_i:= U\cap X_i$. By the definition of $T_X$, we have that $T_X(U_i)=\rho(\ga_i) (U_i)$. But all the maps $\rho(\ga_i)$ are measure-preserving, and so the claim follows.
\end{proof}

Let $(X,T_X)$ be a Turing dynamical system. If $x\in X$ is such that for some $k$ we have $T_X^k(x)= T_X^{k+1}(x)$ then put $T_X^\infty(x):= T_X^k(x)$. Otherwise leave $T_X^\infty$ undefined.

The {\bf first fundamental set}, or simply the {\bf fundamental set}  of $(X,T_X)$ is the subset $\call F_1(T_X)$ of $I$ consisting of all those points $x$ such that $T^\infty(x) \in A$ and for no point $y\in X$ one has $T_X(y)=x$.  The {\bf second fundamental set} of $(X,T_X)$ is the subset $\call F_2(T_X):= T^\infty (\call F_1(T_X))$.  Both the first and the second fundamental set of $(X,T_X)$ are measurable. Indeed, for example, the first fundamental set is nothing but  
$$
	\bigcup_{i=1}^\infty (T_X^{-i}(A)\cap I)-T_X(X),
$$
and $T_X(X)$ is measurable by Proposition \ref{prop_measure_contracting}. Therefore, we define the {\bf first fundamental value}, or simply the {\bf fundamental value}, of $(X,T_X)$ as the measure of its first fundamental set, and similarly for the {\bf second fundamental value}. They are denoted by $\Om_1(T_X)$ and $\Om_2(T_X)$.

We  say that $(X,T_X)$ {\bf stops on any configuration}, if for almost any $x$ $T_X^\infty(x)\ A\cup R$; it {\bf has disjoint accepting chains}, if for almost all different points $x, y\in \call F_1(T_X)$ we have $T^\infty(x) \neq T^\infty(y)$; finally it {\bf does not restart}, if $\mu(T_X(X)\cap I)=0$. 

\begin{proposition}\label{prop_equality_of_fundamental_values}
If $(X,T_X)$ has disjoint accepting chains, then $\Om_1 = \Om_2$.
\end{proposition}
\begin{proof}
By assumption, the map $T_X^\infty\colon\call  F_1 \to \call F_2$ is injective almost everywhere. We can express $\call F_1$ as a countable disjoint union
$$
	\call F_1 = \bigcup_{i=1}^\infty T_X^{-i}(\call F_2)\cap \call F_1,
$$
and it is clear that $T_X^\infty$ restricted to $T_X^{-i}(\call F_2)\cap \call F_1$ is equal to $T_X^i$. Therefore, the claim follows by Proposition \ref{prop_measure_contracting}.
\end{proof}

\subsec{Expressing the fundamental values as von Neumann
dimensions}\label{subsec_expressing_the_fundamental_values}
 
Let $(X,T_X)$ be a Turing dynamical system. Let $\call G$ denote the groupoid $\call G(\rho)$ and define $T$ to be an element of the groupoid ring $\C \call G$ given by $T:= \sum_{i=1}^n \ga_i\chi_i$, where $\chi_i$'s are characteristic functions of the respective $X_i$'s. 

Given two operators $A$ and $B$ on a Hilbert space, note that $\ker A^*A + B^*B =\ker A \cap \ker B$.

Let $S\in \C \call G$ be defined as 
$$
S:=  (T+ \chi_X - \chi_I - \chi_A - \chi_R)^*(T+ \chi_X - \chi_I - \chi_A - \chi_R) + \chi_A.
$$

\begin{theorem}\label{main_theorem}
If $(X, T_X)$ stops on any configuration and doesn't restart then $\dimvn\ker S$ is equal to $\mu(I) - \Om_2(T_X)$.
\end{theorem}

\begin{proof}
Let $\Phi$ be the sequence ${\ga_i}_{|X_i}$ of measurable edges. In order to use Proposition \ref{prop_computational_tool} we show that $\call G(\Phi)$ is a relation groupoid with finite orbits.

\begin{lemma*} (1) The orbit $\call G(\Phi)x$ of almost any point $x$ is finite. (2) For almost every point $x$ the $\Phi$-diagram $\schr_\Phi(x, \call G(\Phi))$ is a tree with one self-loop labeled by the neutral element of $\Ga$ restricted to some set.
\end{lemma*}
\begin{proof}
Note that in $\schr_\Phi(x, \call G(\Phi))$ there is an oriented edge between two points $y$ and $z$ precisely when $T_X(y) = z$. Therefore the second statement follows. 

\begin{figure}[h]%
  \resizebox{0.8\textwidth}{!}{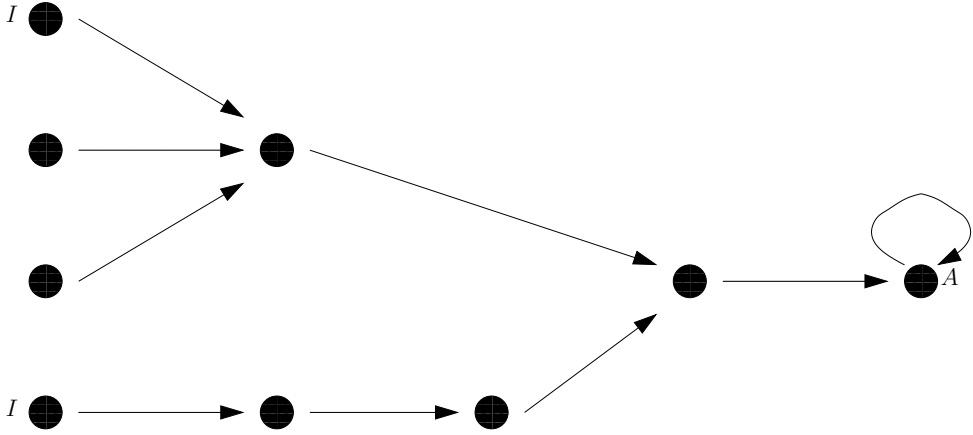}
  \caption{An example Schreier graph of a point. Some of the starting points belong to $I$, and the final point belongs to either $A$ or $R$. The only loop in the graph is the self-loop at the final point whose label is the neutral element of $\Ga$.}
  \label{tgraph_example}
\end{figure}

Let $X^0\subset X$ be the set of fixed points of $T_X$. Since the system is strongly attracting, $X^0$ is equal to $A\cup R\cup Z$, where $Z$ is some set of measure $0$.

Similarly define $X^i, i>0$, as 
$$
	X^i := \{ x\in X: \quad T_X^i(x)\in X^0,  T^j(x)\notin X^0 \text{ for } j<i\}.
$$
Since $X^{i+1} = T_X^{-1}(X^{i})-\cup_{j\le i} X^j$, we see inductively that the sets $X^i$ are measurable. It is clear that they are disjoint, and because the system is strongly attracting we know that
$$
	\mu \left(\bigcup_{i\ge 0} X^i\right) = 1
$$ 
In particular $\lim_{j\to +\infty} \mu (X^j) = 0.$, and so the set 
$$
	\bigcap_{j\ge 0} T_X^{j} (X^j).
$$
is of measure $0$. But every infinite $\call G(\Phi)$-orbit intersects the above set. The set of points with infinite $\call G(\Phi)$-orbits is contained in countably many translates of the above set, and so it is also of measure $0$.
\end{proof}

By assumption $A\cup R$ is a fundamental domain for $\call G(\Phi)$ and so by Proposition \ref{prop_computational_tool} we have that
$$
  \dimvn\ker S = \int_{A\cup R} \dim \ker S_x \,d\mu(x)
$$

\begin{lemma*}
Let $x\in X$. Statement of the theorem follows if we show that $\dim \ker S_x$ is equal to $0$ if $|\call G(\Phi)x \cap I|=0$ and 
\begin{equation}\label{eq_formula_I-A}
	|\call G(\Phi)x \cap I | -| \call G(\Phi)x \cap A|
\end{equation}
otherwise.
\end{lemma*}
\begin{proof}
We get
\begin{eqnarray*}
  \dimvn\ker S &=& \int_{A\cup R} \dim \ker S_x \,d\mu(x) \\
	       &=& \int_{T_X^\infty(I)} |\call G(\Phi)x\cap I| - | \call G(\Phi)x \cap A| \,d\mu(x) \\
	       &=& \int_{T_X^\infty(I)} |\call G(\Phi)x\cap I| \,d\mu(x)  -  \int_{T_X^\infty(I)\cap A}  | \call G(\Phi)x \cap A|  \,d\mu(x) \\
	       &=& \int_{T_X^\infty(I)} |\call G(\Phi)x\cap I| \,d\mu(x)  -  \Om_2(T_X), 
\end{eqnarray*}
and the statement follows from Lemma \ref{lemma_integration_change_A_and_B}.
\end{proof}

Let us first compute the kernel of $T_{x}+ (\chi_X)_x - (\chi_I)_x - (\chi_A)_x - (\chi_R)_x$. From the definition $T_x$ is the convolution operator on the Schreier graph $\schr_\Phi(x, \call G(\Phi))$, and the operator $(\chi_X)_x - (\chi_I)_x - (\chi_A)_x - (\chi_R)_x$ acts identically on vectors $\zeta_v$ for $v\notin A\cup I \cup R$ and is null on other vectors $\zeta_v$. 

Let $x_1$ be any starting vertex of $\schr_\Phi(x, \call G(\Phi))$ which belongs to the initial set $I$, and let $(x_1, x_2, \ldots, x_k)$ be such that $T_X(x_i) =x_{i+1}$ when $i<k$ and $T_X(x_k)=x_k$. Computation shows that the vector
$$
  \xi(x_1) := \zeta_{x_1} - \zeta_{x_2} + \ldots  \pm \zeta_{x_k}
$$
is in the kernel of $T_{x}+ (\chi_X)_x - (\chi_I)_x - (\chi_A)_x - (\chi_R)_x$. Furthermore, if $u,v,\ldots,w$ are different starting vertices then the vectors $\xi(u), \xi(v), \ldots, \xi(w)$ are linearly independent. 

\begin{lemma*}
The linear span of the set $\{\xi(x): x\in V(g)\cap I  \}$ is equal to the kernel of $T_{x}+ (\chi_X)_x - (\chi_I)_x - (\chi_A)_x - (\chi_R)_x$
\end{lemma*}
\begin{proof}
We just saw containment. On the other hand, suppose there is a vector $\eta$ in the kernel such that $\langle \eta, \zeta_v \rangle=0$ for every $v\in I$. Let $w$ be a maximal (with respect to the relation $w\ge v \iff T_X^k(w)=v$ for some $k$) vertex of $\call G(\Phi)x$ such that $\langle \eta, \zeta_w \rangle \neq 0$.

Write $\eta = \zeta_w +\eta'$, where $\eta'$ is a linear combination of vectors $\zeta_u$ for $w>u$. It follows that $\langle (T_x+ (\chi_X)_x - (\chi_I)_x - (\chi_A)_x - (\chi_R)_x) \eta', \zeta_w \rangle = 0$, and so
\begin{multline*}
0 = \langle (T_x+ (\chi_X)_x - (\chi_I)_x - (\chi_A)_x - (\chi_R)_x) \eta, \zeta_w \rangle = \\
= \langle (T_x+ (\chi_X)_x - (\chi_I)_x - (\chi_A)_x - (\chi_R)_x) \zeta_w, \zeta_w \rangle. 
\end{multline*} 
But if $w$ is not a final vertex then $\langle T_x \zeta_w, \zeta_w \rangle =0$, and  $\langle (\chi_X)_x - (\chi_I)_x - (\chi_A)_x - (\chi_R)_x) \zeta_w, \zeta_w\rangle  = 1$ (the last equality follows, since by assumption also $w\notin I$), which is a contradiction. And if $w$ is a final vertex then $\langle T_x \zeta_w, \zeta_w \rangle =1$, and $\langle ((\chi_X)_x - (\chi_I)_x - (\chi_A)_x - (\chi_R)_x) \zeta_w, \zeta_w \rangle = 0$, which  also is  a contradiction. 
\end{proof}

Now we need to consider two cases: the final point of $\schr_\Phi(x, \call G(\Phi))$ is in $R$ or $A$. If it is in $R$ then $(\chi_A)_x=0$ and $\ker S_x = \ker  T_{x}+(\chi_X)_x - (\chi_I)_x - (\chi_A)_x - (\chi_R)_x$. Therefore the formula \ref{eq_formula_I-A} holds because the dimension of linear span of $\{\xi(v): v\in I  \}$ is precisely $|I|$.  

If the final point of $\schr_\Phi(x, \call G(\Phi))$ is in $A$ then kernel of $(\chi_A)_x$ is of codimension $1$, and it non-trivially intersects span of $\{\xi(x): x\in I  \}$, as soon as the latter set is non-empy. This shows that the formula \eqref{eq_formula_I-A} holds also in this case.
\end{proof}

\begin{corollary} \label{main_cory} If $(X, T_X)$ stops on any configuration, is strongly attracting, and has disjoint accepting chains, then $ \dim_{vN} \ker S = \mu(I) - \Om_1(X,T_X).$ 
\end{corollary}
\begin{proof}
 Follows from Proposition \ref{prop_equality_of_fundamental_values}.
\end{proof}

\section{Turing dynamical systems - examples}\label{section_tds-examples}

We want to use presented examples later together with Theorem \ref{main_theorem} and Proposition \ref{prop_pontryagin}, so we have to assure that the actions are by continuous group automorphisms, which is the reason for an extra degree of complicacy.

\subsec{Turing dynamical system associated to a set of natural numbers}\label{subsec_first_example} 

\paragraph{Definition of $X$ and $\Ga$} Let $X$ be a measure space $M^\Z \times S$, where $M := \Zmod{2}\oplus\Zmod{2}\oplus\Zmod{2}$ should be
interpreted as the set of symbols, and $S:=\Zmod{2}\oplus\Zmod{2}\oplus\Zmod{2}$ as the set
of  states of a Turing machine.

Let $\Aut(M)$ be the group of group automorphisms of $M$; similarly for $\Aut(S)$. Recall that $\Aut(M) \wr \Z$ is defined as $\Aut(M)^{\oplus \Z} \rtimes \Z$, and put 
$$
\Ga = [(\Aut(M)\wr\Z) \ast \Zmod{2}]\times \Aut(S),
$$
where $\ast$ denotes the free product.

\paragraph{Notation for elements of $X$} Let $m_{-k} \ldots m_0 \ldots m_l \in M$ and $\si \in S$. The set
$$
\{(\dot n_i, \tau) \in M^\Z \times S: n_{-k} = m_{-k}, \ldots, n_l = m_l, \,\, \tau = \si\}
$$ 
is denoted by
$$
[m_k\, m_{k-1} \ldots  m_{-1}\, \underline{m_0}\, m_1  \ldots m_l][\si].
$$ 
Given $\si\in S$, let 
$$
[][\si]:= \bigcup_{m\im M} [\underline m][\si],
$$
 and given $m\in M$, let 
$$
[\underline m][]:= \bigcup_{\si\in S} [\underline m][\si].
$$ 

A concrete element from the set $[m_k\, m_{k-1}\ldots  m_{-1}\, \underline{m_0}\, m_1  \ldots m_l][\si]$ is denoted by 
$$
(m_k\, m_{k-1} \ldots  m_{-1}\, \underline{m_0}\, m_1  \ldots m_l)[\si],
$$
and similarly for $()[\si]$ and $(\underline m)[]$.

\paragraph{Definition of the action} Let us fix a set of positive natural numbers $\Si$. The action $\rho$ of $\Ga$ on $X$ will depend on $\Si$. Whenever we want to stress this dependence we use the symbol $\rho_\Si$.

Let $\be'$ be the unique automorphism of $M$ such that  $\be'((1,0,0)) = (1,0,0)$, $\be'((0,1,0)) = (0,0,1)$, $\be'((0,0,1))=(0,1,0)$. Let $M^\be$ be the set of fixed points of $\be'$. Note that it consists of $4$ elements. Let $\beta$ be the automorphism of $M^\Z$ defined by
\begin{equation}
	\big(\beta(\dot x_i)\big)_j := \left\{ 
\begin{array}{l l}
\be'(x_0)  & \quad \mbox{if $j= 0$ }\\
  x_j  & \quad \mbox{otherwise}\\
\end{array}\right.
\end{equation}	
The automorphism $\be$ will be referred to as {\bf normal flip}. 

Let $B=B_\Si$  be the following automorphism of $M^\Z$:
\begin{equation}
	\big(B(\dot x_i)\big)_j := \left\{ 
\begin{array}{l l}
x_j  & \quad \mbox{if $j\notin \Si$ }\\
  \be'(x_j)  & \quad \mbox{otherwise}\\
\end{array}\right.
\end{equation}	
The automorphism $B$ will be referred to as {\bf oracle flip}.

We proceed to describe $\rho$. The subgroup $\Aut(M)\wr \Z <\Ga $ acts in the standard way on the $M^\Z$ coordinate of $X$: $\Aut(M)$ acts on the $0$-coordinate of $M^\Z$ in the natural way, and the generator $t$ of $\Z$ acts by 
\begin{equation}
	\big(\rho(t)(\dot m_i)\big)_j := m_{j+1}.
\end{equation}	
The maps  $\rho(t)$ and $\rho(t^{-1})$ will be called respectively {\bf shift forward}, and  {\bf shift backward}.

The subgroup $\Aut(S) < \Ga$ acts on the $S$ coordinate of $X$ in the natural way.

The generator of $\Zmod{2}$ acts by the oracle flip $B$ on $M^\Z$; it will be also denoted by $B$.

\paragraph{Division of $X$ into disjoint subsets} We choose the following division: 
$$
X = \bigsqcup_{m\in M, \si \in S} [\underline m][\si].
$$

We just finished defining a dynamical hardware $(X)$. When we need to stress its dependence on $\Si$, we denote it by $(X^\Si)$.

Elements of $S = \Zmod{2}\oplus\Zmod{2}\oplus\Zmod{2}$ will be referred to as [\textit{Start}],
[\textit{Search forward for $(0,1,0)$}], [\textit{Search backward for $(0,1,0)$}], [\textit{Search forward
for either $(0,1,0)$ or $(0,0,1)$}], [\textit{Dummy state 1}], [\textit{Dummy state 2}], [\textit{Dummy state 3}],
[\textit{Dummy state 4}]. We do not specify which names correspond to which elements of $S$ - the only
important thing is that this assignment is made in such a way that $\Aut(S)$ acts transitively on the first four elements.

We proceed to define a dynamical software for $(X^\Si)$. It will be denoted by $(T_X)$ or $(T^\Si_X)$.

\begin{figure}[h]%
  \resizebox{0.8\textwidth}{!}{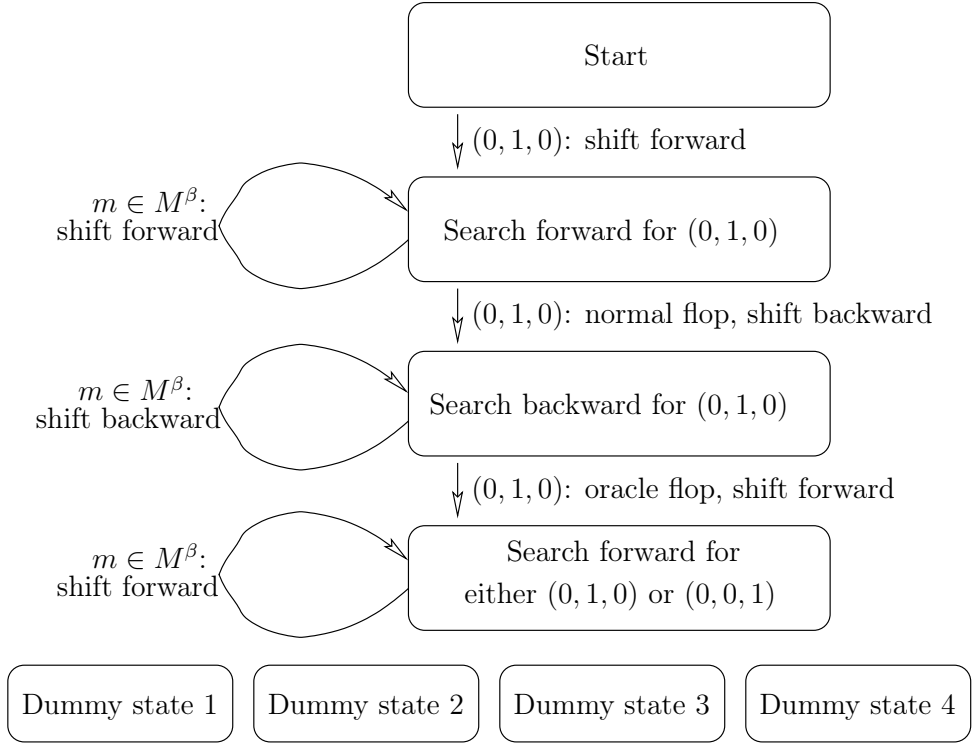}
  \caption{Turing dynamical system $(X^\Si, T^\Si_X)$}
  \label{fig_algorithm}
\end{figure}

\paragraph{Choice of elements of $\Ga$ for the sets $[\underline m][\si]$} This is done in Figure \ref{fig_algorithm}: arrow between two
states $\si$ and $\tau$ with a label, for example, 
\begin{center}
$(0,1,0)$: normal flip, shift backward
\end{center}
means that the element of $\Ga$ corresponding to $[\underline{(0,1,0)}][\si]$ is $a \cdot \be \cdot t^{-1}$,
where $a\in \Aut(S)$ is any group automorphism sending $\si$ to $\tau$. Since the action is a right action it means in particular that in this example flipping is done before shifting.
Similarly, an arrow with label
\begin{center}
 $m\in M^\be$: shift forward
\end{center}
joining a state $\si$ with itself means that the element of $\Ga$ corresponding to $[\underline{m}][\si]$ for $m\in
M^\be$ is $t$.

Finally, if for some state $\si$ there is no label with a given symbol $m\in M$ then it means that the element
of $\Ga$ corresponding to $[\underline{m}][\si]$ is the neutral element. 

\paragraph{Choice of the sets $I$, $A$ and $R$}. We specify them as follows:
\begin{eqnarray*}
	I &:=& \big[\underline{(0,1,0)}\big]\big[\textit{Start}\big],\\
	A &:=& \big[\underline{(0,1,0)}\big]\big[\textit{Search forward for $(0,1,0)$ or $(0,0,1)$}\big].
\end{eqnarray*}
As to the set $R$, it is defined as the union of all the sets $[\underline m][\si]$ whose associated group
element is the neutral element, apart from $A$.

\subsec{Properties of the system $(X,T_X)$}\label{subsec_first_properties}

\begin{proposition}\label{prop_fundamental_set}
 The first fundamental set of $(X, T_X)$ is equal to  the union
 $$
 		\bigcup_{k\in \Si} F_k \cup Z,
 $$
where $F_k$ is  equal to
$$
	\bigcup_{m_1,\ldots ,m_{k-1}\in M^\be}
\big[\underline{(0,1,0)}\,m_1\,m_2\ldots m_{k-1}\,(0,1,0)\!\big]\big[\!\textit{Start}\big],
$$
and $Z$ is some set of measure $0$.
\end{proposition}

\begin{proof} This proposition follows from chasing through Figure \ref{fig_algorithm}. First, we show that $F_k$ is in the fundamental set for $k\in\Si$.

Let $x=(\underline{(0,1,0)}\, m_1 \, m_2 \ldots m_{k-1} \,(0,1,0))[\textit{Start}]$, $m_i\in M^\be$, $k\in \Si$.
Because of the arrow between the first and the second level of Figure 1, we have 
$$
T(x) = \big(\!{(0,1,0)}\, \underline{m_1} \, m_2 \ldots m_{k-1} \,(0,1,0)\!\big)\big[\!\textit{Search forward for $(0,1,0)$}\!\big].
$$
Then, because of the arrow ``$m \in M^\be\colon$ shift forward'' on the second level of Figure 2, we get
$$
T^{k}(x) = \big(\!{(0,1,0)}\, m_1 \, m_2 \ldots {m_{k-1}} \,\underline{(0,1,0)}\big)\big[\!\textit{Search forward for $(0,1,0)$}\!\big].
$$
Because of the arrow between the second and the third level, we see
$$
T^{k+1}(x) = \big(\!{(0,1,0)}\, {m_1} \, m_2 \ldots \underline{m_{k-1}} \,(0,0,1)\!\big) \big[\!\textit{Search backward for $(0,1,0)$}\!\big].
$$
Because of the arrow ``$m \in M^\be\colon$ shift backward'' on the third level, we conclude
$$
T^{2k}(x) = \big(\underline{(0,1,0)}\, m_1 \, m_2 \ldots m_{k-1} \,(0,0,1)\!\big) \big[\!\textit{Search backward for $(0,1,0)$}\!\big].
$$
Because of the arrow between the third and the fourth level, and since $k\in \Si$, we get that $T^{2k+1}(x)$ is equal to
$$
\big(\!(0,1,0)\, \underline{m_1} \, m_2 \ldots m_{k-1} \,(0,0,1)\!\big) \big[\!\textit{Search forward
for either $(0,1,0)$ or $(0,0,1)$}\!\big].
$$
Finally, because of the arrow ``$m \in M^\be\colon$ shift forward'' on the fourth level, we see that $T^{3k}(x)$ equals
$$
\big(\!(0,1,0)\, {m_1} \, m_2 \ldots m_{k-1} \,\underline{(0,0,1)}\big) \big[\!\textit{Search forward
for either $(0,1,0)$ or $(0,0,1)$}\!\big],
$$
which is an element of the accepting set.

\medskip

In the other direction, let
$$
	x = (\underline{m_0}\,m_1\ldots m_{k-1}\,m_{k})[\si]
$$
be an element of the fundamental set of $(X,T_X)$. We can assume $m_1,\ldots, m_{k-1} \in M^\be$, and
$m_k\notin M^\be$, since points $(\dot x_i,\si)\in X$ such that $x_i\in M^\be$ for $i>0$ form a set of
measure $0$. We have to prove that (1) $\si = $[\textit{Start}], (2) $m_0= (0,1,0)$, (3) $m_k=(0,1,0)$,
and (4) $k\in \Si$.

(1) and (2)  follow from $I=\underline{(0,1,0)}\,[\textit{Start}]$.  As before we have 
$$
	T_X^k(x) = ({m_0}\,m_1\,\ldots\, {m_{k-1}}\,\underline{m_{k}})[\textit{Search forward for $(0,1,0)$}],
$$
and therefore from the fact that $T_X^{k}(x) \notin R$ we get (3). Again, as before we see 
$T_X^{2k}(x)$ is equal to
\begin{equation*}\label{crucial_moment}
	\big(\underline{m_0}\,{m_1}\,\ldots\, {m_{k-1}}\,(0,0,1)\!\big)\big[\!\textit{Search backward for $(0,1,0)$}\!\big].
\end{equation*}
Now, suppose that (4) does not hold, i.e.  $k\notin \Si$. Then because of the arrow between the third and the fourth level, and by the definition of the oracle flip, we get that $T_X^{2k+1}(x)$ is 
$$
	\big(\!{m_0}\,\underline {m_1}\,\ldots\, {m_{k-1}}\,(0,0,1)\!\big) \big[\!\textit{Search forward for either $(0,1,0)$ or $(0,0,1)$}\!\big],
$$
which implies that $T_X^{3k}(x)$ is 
$$
	\big(\!{m_0}\,{m_1}\,\ldots\, {m_{k-1}}\,\underline{(0,0,1)}\big)\big[\!\textit{Search forward for either $(0,1,0)$ or $(0,0,1)$}\!\big],
$$
which is an element of $R$, which contradicts the assumption that $x$ is in the fundamental set.
\end{proof}

\begin{cory}\label{cory_fundamental_value}
The first fundamental value of $(X, T_X)$ is equal to
$$
	\frac2{8^3} \sum_{i\in \Si} \frac1{2^i}
$$
\end{cory}

\begin{proof}
The sets $F_k$ are  disjoint, and the measure of $F_k$ is equal to $\frac18
(\frac12)^{k-1}\frac18\frac18 = \frac2{8^3} \frac1{2^k}$.

\end{proof}

\begin{proposition}\label{prop_tds_properties}  The Turing system $(X, T_X)$ (i) doesn't restart, (ii) has disjoint accepting chains, and (iii) stops on any configuration.
\end{proposition}
 
\begin{proof}
(i) We have to show that $T_X(I) \cap I = \emptyset$.  Recall $I =[\underline{(0,1,0)}][\textit{Start}]$; from Figure \ref{fig_algorithm} we see that (1) points from outside of $[][\textit{Start}]$ are not mapped into $[][\text{Start}]$, and in particular are not mapped into $I$; (2) points from $I$ are mapped outside of $[][\textit{Start}]$ and in particular outside of $I$; and (3) Points from $[][\textit{Start}]$ which are not in $I$ are mapped identically to themselves, and so are also not mapped to $I$. 

(ii) We need to check that there exists a subset of the fundamental set, which has the same measure as the fundamental set, and on which the map $T_X^\infty$ is injective. But in the proof of Proposition \ref{prop_fundamental_set} we saw that if $x\in F_k$, $k\in \Si$, and we write $x = (\underline{(0,1,0)}\, m_1 \, m_2 \ldots m_{k-1} \,(0,1,0))[\textit{Start}]$ then $T_X^\infty(x)$ is equal to $T^{3k}(x)$, which is
$$
((0,1,0)\, {m_1} \, m_2 \ldots m_{k-1} \,\underline{(0,0,1)}) [\textit{Search forward for either $(0,1,0)$ or $(0,0,1)$}],
$$
from which we can recover $x$. In particular $T_X^\infty$ is injective on $\bigcup_{k\in \Si} F_k$.

For (iii), note the following claim.
\begin{claim*}
 Suppose $x= (\dot m_i, \si) \in M^\Z \times S = X$ is such that $T^\infty_X(x) \notin A\cup R$ or is undefined. Then there exists $N\in \Z$ such that $m_j\in M^\be$ for all $j\ge N$ or for all $j\le N$.
\end{claim*}
\begin{proof}
The only elements of $\Ga$ which are applied to the $M^\Z$-coordinate are shifts, normal flip and oracle flip. Since they preserve the property ``exists $N$ such that $m_j\in M^\be$ for all $j \ge N$ or all $j \le N$'', it is enough to show that some power $T_X^j(x)$ has this property. But this is clear from Figure 1: first, there exists a state $\tau$ and $K$ such that $T_X^k(x)\in [][\tau]$ for $k \ge K$. But the only way it is possible is if either (1)  for some $k$, $T_X^k(x)$ is in a set $[m][\si]$ for which the corresponding element of $\ga$ is the neutral element, or (2)  if we write  $T^K(x) = (\dot n_i, \tau)$ then exists $N$ such that $n_j\in M^\be$ for all $j> N$ or all $j<N$. The second case is what we want to show, and the first case is not possible because it implies that $T^{k}(x) \in A\cup R$.
\end{proof}
Now, (iii) follows, since the set of points $\dot m_i\in M^\Z$ such that there exists $N$ such that $m_j\in M^\be$ for $j> N$ is of measure $0$.
\end{proof}

\subsec{A ``read only'' system with irrational fundamental values}

\label{subsec_second_example}

\paragraph{Definition of $Y$ and $\De$} Let $Y$ be the measure space $(\Zmod{2}^\Z$ $\times$ $\Zmod{2}^\Z$ $\times$ $\Zmod{2}^\Z)${}$\times${}$S$, where
$S=\Zmod{2}\times\Zmod{2}\times\Zmod{2}$ should again be interpreted as the set of states. However, in this example $\Zmod{2}$ should be interpreted as the set of symbols, and $\Zmod{2}^\Z \times \Zmod{2}^\Z \times \Zmod{2}^\Z$ should be interpreted as the set of \textit{triples} of tapes.

Let $\De$ be the group $\Z^3\times \Aut(S)$, and let standard generators of $\Z^3$ be denoted by $t_1$, $t_2$, and
$t_3$. They will be also referred to as {\bf shift forward tape 1}, {\bf shift forward
tape 2 }, and {\bf shift forward tape 3}. Similarly $t_i^{-1}$ will be referred to as {\bf shift backward tape
$i$}.

\paragraph{Definition of the action} The action $\rho: \De \curvearrowright Y$ is defined as follows. Let $V: \Zmod{2}^\Z \to \Zmod{2}^\Z$ be the shift automorphism, i.e. $(V(\dot x_i))_j = x_{j+1}$. Define
\begin{eqnarray*}
	\rho(t_1)(\dot x_i, \dot y_i, \dot z_i) &:=& (V(\dot x_i), \dot y_i, \dot z_i) \\
	\rho(t_2)(\dot x_i, \dot y_i, \dot z_i) &:=& (\dot x_i, V(\dot y_i), \dot z_i) \\
	\rho(t_3)(\dot x_i, \dot y_i, \dot z_i) &:=& (\dot x_i, \dot y_i, V(\dot z_i)).
\end{eqnarray*}
$\Aut(S)$ acts in the natural way on the $S$-coordinate.

\paragraph{Notation} Notation is similar to that in Subsection \ref{section_tds-examples}\ref{subsec_first_example}. Given $a_{-k_a}, \ldots, a_{l_a}$, $b_{-k_b}, \ldots, b_{l_b}$, $c_{-k_c},\ldots  c_{l_c} \in \Zmod{2}$ and $\si \in S$, the set
\begin{multline*}
\{(\dot x_i,\dot y_i,\dot z_i,\tau) \in (\Zmod{2}^\Z \times \Zmod{2}^\Z \times \Zmod{2}^\Z) \times S \colon \\ x_{-k_a} = a_{-k_a} ,\ldots, z_{l_c} = c_{l_c}, \tau = \si\}
\end{multline*}
is denoted by
$$
\left[ \begin{array}{lll}
 	a_{-k_a}\,\ldots\, &\underline{a_0}&\, \ldots \, a_{l_a}\\	
 	b_{-k_b}\,\ldots\, &\underline{b_0}&\, \ldots \, b_{l_b}\\	
 	c_{-k_c}\,\ldots\, &\underline{c_0}&\, \ldots \, c_{l_c}\\	
	\end{array}
\right] [\si].
$$

Given $v = (v_1, v_2, v_3) \in \Zmod{2}^3$ and $\si\in S$, the set
$$
\left[ \begin{array}{lll}
 	\underline{v_1} \\	
 	\underline{v_2}\\	
 	\underline{v_3}\\	
	\end{array}
\right] [\si].
$$
is also denoted by $[\underline v][\si]$ and $[\underline{(v_1,v_2,v_3)}][\si]$.

Within the context of above notation, given a natural number $k$ and $\eps\in \Zmod{2}$ we denote by $\eps^k$ the sequence of $k$ consecutive $\eps$'s. 

\paragraph{Division of $Y$} We choose the following decomposition of $Y$: 
$$
	Y=\bigsqcup_{v\in\Zmod{2}^3,\, \si\in S} [\underline{v}][\si].
$$ 

We have just defined a dynamical hardware $(Y)$

The states - i.e. elements of $S$ - will be referred to as [\textit{Start}],
[\textit{Check if the number of $0$'s on tape 1 and on tape 2 is the same}], [\textit{Search backward for $1$
on tape $1$}], [\textit{Search forward for $1$ on tape 1 and move forward tape 3}], [\textit{Dummy state 1}],
[\textit{Dummy state 2}], [\textit{Dummy state 3}], [\textit{Dummy state 4}]. Again we do not specify which of
these names correspond to which elements of $S$, but we demand that $\Aut(S)$ acts transitively on the first
four of the above states.

\begin{figure}[h]%
  \resizebox{0.9\textwidth}{!}{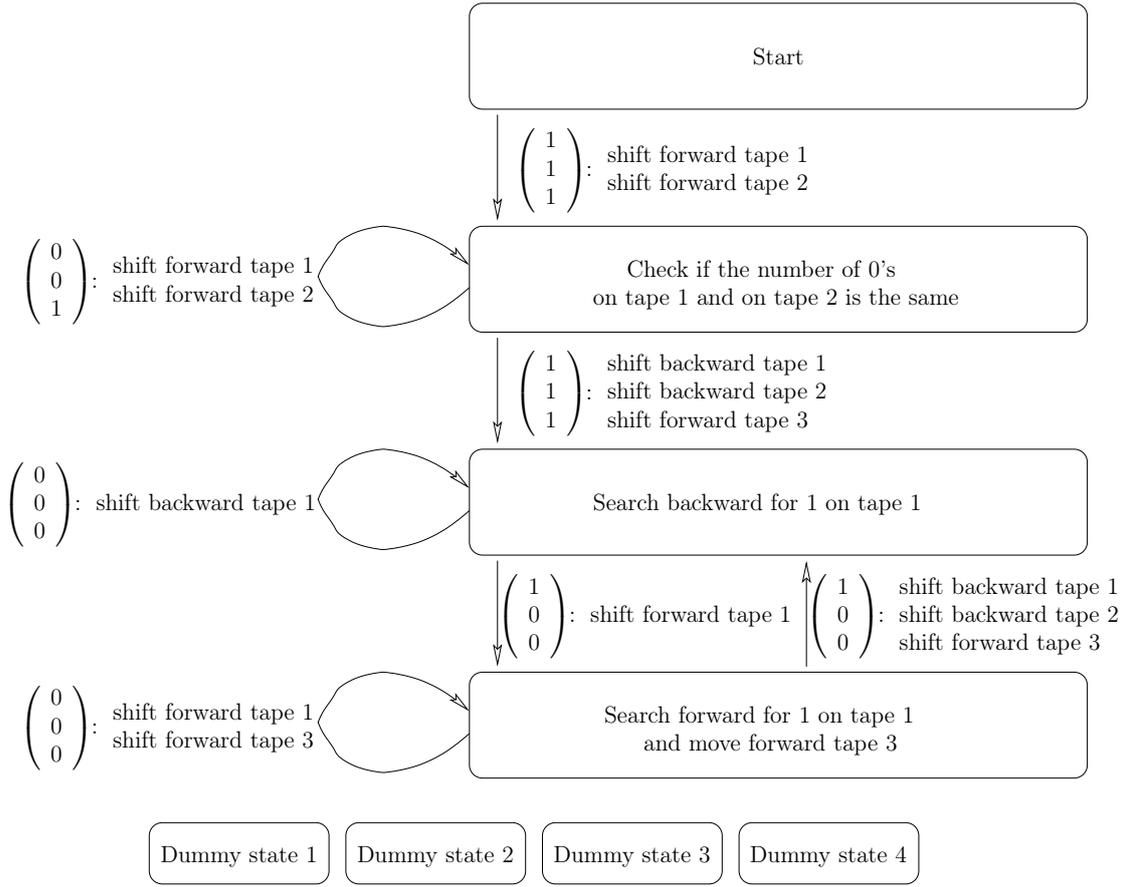}
  \caption{Turing dynamical system $(Y, T_Y)$}
  \label{fig_algorithm2}
\end{figure}

We proceed to define a dynamical software $(T_Y)$ for $(Y)$. This is done using Figure 2, with the
same convention as in Subsection \ref{section_groupoids_examples}\ref{subsec_first_example}, e.g. an arrow with a
label
$$
\left( \begin{array}{c}
	1 \\
	1 \\
	1 \end{array}\right)
\!\!:\! 
\left.\begin{array}{ll} 
	\text{shift forward tape 1} \\
	\text{shift forward tape 2}
\end{array}\right.
$$  
between states $\si$ and $\tau$ means that the element of $\De$ associated to the set $[\underline{(1,1,1)}][\si]$ is
$a\cdot t_1\cdot t_2$, where $a\in \Aut(S) $ is any group automorphism of $S$ which sends $\si$ to $\tau$. 

\paragraph{Choice of the sets $A$, $I$ and $R$} We define them as follows:
\begin{eqnarray*}
	I &:=& [\underline{(1,1,1)}][\textit{Start}],\\
	A &:=& [\underline{(0,1,1)}][\textit{Search backward for $1$ on tape 1}],
\end{eqnarray*}
and the rejecting set $R$ is defined to be the union of all the sets $[\underline{v}][\si]$  whose associated group
element is the neutral element, apart from $A$.

\medskip

\subsec{Properties of the system $(Y,T_Y)$}

\begin{proposition}\label{prop_2fundamental_set}
 The first fundamental set of $(Y, T_Y)$ is equal to  the union
 $$
 		\bigcup_{k=1}^\infty  F_k \cup Z,
 $$
where $Z$ is some set of measure $0$ and
$$
F_k:=\left[ \begin{array}{l}
 	\underline{1}\,0^k \, 1\\	
 	\underline{1}\,0^k \, 1\\
 	\underline{1}\,0^{k^2+2k}\, 1
	\end{array}
\right] [\textit{Start}].
$$
\end{proposition}

\begin{proof} The proof is by chasing Figure \ref{fig_algorithm2}, fully analogous to the proof of Proposition \ref{prop_fundamental_set}.
\end{proof}

\begin{cory}\label{cory_2fundamental_value} The first fundamental value of $(Y, T_Y)$ is equal to
 $$
	\frac18\sum_{k=1}^\infty \frac1{2^{k^2+4k+6}}
$$
\end{cory}
\begin{proof}
 Indeed, measure of the set $F_k$ is equal to $\frac18\cdot \frac1{2^{k^2+2k}} \cdot \frac1{2^k}\cdot
\frac1{2^k} \cdot\frac1{2^6}$ (i.e. ``$\frac18$ is for states, $\frac1{2^{k^2+2k}}$ and $\frac1{2^k}$ are for
$0$'s, and $\frac1{2^6}$ is for $1$'s''.)
\end{proof}

\begin{proposition}\label{2oracle_reductions}  The Turing system $(Y, T_Y)$ (i) doesn't restart, (ii) has disjoint accepting chains, and (iii) stops on any configuration.
\end{proposition}

\begin{proof}
(i) and (ii) are proved just like (i) and (ii) of Proposition \ref{prop_tds_properties}. 

As for (iii), let $y\in Y$ be a point such that $T_Y^k(y)\notin A\cup R$, and let $y_i := T_Y^k(y)$. Let $\{y_i\}_{i\ge k}$ denote the set of elements of the sequence $y_k, y_{k+1}, \ldots$, and $\{y_i\}:= \{y_i\}_{i\ge 0}$. 

If $\{y_i\}$ is contained in the first two levels of Figure \ref{fig_algorithm2} then there are infinitely many consecutive $0$'s to the right on the first two tapes of $y$. The set of points which have infinitely many consecutive $0$'s on any tape has measure $0$; and if $y_k$ is on the third or fourth level then $\{y_i\}_{i\ge k}$ is contained in the third and the fourth level. Therefore swe can assume that $\{y_i\}_{i\ge l}$ is for some $l$ contained in the third and the fourth level of Figure \ref{fig_algorithm2}.

Consider two possibilities: (1) $\{y_i\}_{i>k}$ is contained in the third level for some $k$, or (2) $y_i$ is in the fourth level for infinitely many $i$'s.

If (1) holds then because of the arrow
$$
\left( \begin{array}{c}
	0 \\
	0 \\
	0 \end{array}\right)
\!\!:\! 
\left.\begin{array}{ll} 
	\text{shift backward tape 1} 
\end{array}\right.
$$  
we see that the first tape of $y$ has infinitely many consecutive $0$'s to the left.

If (2) holds then note that the only element of $\De$ which acts on the third tape of $y$ is shift forward. From Figure \ref{fig_algorithm2} we see that shift forward on the third tape is applied infinitely many times. We see also that it can be applied only if there is $0$ on the third tape. Therefore in this case there are infinitely many consecutive $0$'s on the third tape of $y$.
\end{proof}

\section{Atiyah problem}\label{section_Atiyah}

\subsec{Preliminaries}

Let $G$ be a discrete countable group. Recall that a real number $r$ is an { $l^2$-Betti number arising from $G$} if for some $k$ there exists $T\in M_k(\Q G)$ such that 
$$
	\dim_{vN} \ker \theta = r.
$$
The set of $l^2$-Betti numbers arising from $G$ is denoted by $\call C(G)$. Note that if $r\in \call C(G)$ then we can always assume that $r=\dimvn\ker T$ for a positive self-adjoint $T$ of norm $<1$, since $\ker T = \ker T^*T = \ker \frac1N T^*T$ for every positive number $N$.

A set $\Si$ of natural numbers is called {\bf computable} if there exists a Turing machine which lists elements of $\Si$ in the increasing order. Equivalently (see e.g. \cite{Lyndon_Schupp}), $\Si$ is computable if there exists a Turing machine which lists elements of $\Si$ in some order (possibly with repetitions) and there exists a Turing machine which lists elements of the complement of $\Si$ in some order (possibly with repetitions).

We say that a real number $r$ has {\bf computable binary expansion} if the fractional part of $r$ is of the form
$$
\sum_{i\in \Si} \frac{1}{2^i}
$$ 
for some computable set $\Si$.

In this section we prove Theorems \ref{thm_lamplighter}, \ref{thm_finitely_generated}, \ref{thm_finitely_presented}, and Corollary \ref{cory_lamplighter}: 

\begin{theorem*}
The set of $l^2$-Betti numbers arising from the group $(\Zmod{2}\wr \Z)^3$ contains 
$$
\frac1{64}- \frac18\sum_{k=1}^\infty \frac1{2^{k^2+4k+6}}.
$$
\end{theorem*}

\begin{theorem*}
The set of $l^2$-Betti numbers arising from finitely generated groups is equal to the set of non-negative real numbers.
\end{theorem*}

\begin{theorem*}
The set of $l^2$-Betti numbers arising from finitely presented groups contains all numbers with computable binary expansions. 
\end{theorem*}

\begin{cory*}
 Let $G$ be a group given by the presentation 
$$
	\langle a,t,s\,|\,a^2=1,[t,s]=1,[t^{-1}at,a]=1,s^{-1}as=at^{-1}at\rangle.
$$ 
The set of $l^2$-Betti numbers arising from $G^3$ contains 
$$
\frac1{64}- \frac18\sum_{k=1}^\infty \frac1{2^{k^2+4k+6}},
$$
\end{cory*}

Corollary follows from the next lemma and the fact proven in \cite{Grigorchuk_Linnel_Schick_Zuk}namely $\Zmod{2} \wr \Z$ is a subgroup of $G$.

\begin{lemma}\label{lemma_subgroup}
 Let $H$ be a subgroup of a group $G$. Then $\call C (H) \subset \call C(G)$.
\end{lemma}
\begin{proof}
 The map $M_k(\C H) \mapsinto M_k(\C G)$, induced by the inclusion $H\mapsinto G$, is a trace-pre\-serv\-ing ${}^*$-homomorphism. The claim follows from Lemma \ref{lemma_trace_preserving_homo}.
\end{proof}

\begin{lemma}\label{lemma_finite_group}
 If $G$ is a discrete countable group and $H$ is a finite group then $|H|\cdot \call C(G\times H) = \call C (G)$.
\end{lemma}
\begin{proof}
 Let $\pi\in \Q(G\times H)$ be the sum 
$$
\frac{1}{|H|}\sum_{h\in H} h
$$
Clearly $\pi$ is a projection of trace $\frac{1}{|H|}$ which commutes with $\Q G$. Similarly the $k\times k$ matrix $\pi_k$ which has $\pi$ everywhere on the diagonal and $0$'s elsewhere is of trace $\frac{1}{|H|}$ and commutes with $M_k (\Q G)$. We claim that for a positive self-adjoint $T\in M_k( \Q G)$ of norm at most $1$ we have
$$
  \frac{1}{|H|}\cdot \dimvn\ker T = \dimvn \ker (1-\pi_k +\pi_k T)
$$
Indeed, by functional calculus the right side is equal to the limit of 
$$
\tau_{G\times H} \left((1 - (1-\pi_k + \pi_kT))^n\right), 
$$
for $n$ going to $+\infty$. The above expression is equal to $\tau_{G\times H}(\pi_k(1+T)^n) = \tau_H (\pi)\cdot \tau((1+T)^n)$. By functional calculus again limit of the latter expressions is equal to $\frac{1}{|H|} \cdot \dimvn \ker T$. This shows the inclusion $C(G) \subset  |H|\cdot \call C(G\times H)$

For the other inclusion, note that the regular representation of $H$ induces a ${}^*$-embedding $\Q H \mapsinto M_{|H|}(\Q)$ such that for $T\in \Q H$ we have $|H|\tau_H(T) = \tr (T)$. This induces a ${}^*$-embedding 
$$
\iota\colon \Q (G\times H) \cong \Q(G) \otimes \Q(H) \mapsinto \Q(G) \otimes M_{|H|}(\Q) \cong M_{|H|}(\Q(G))
$$ 
such that for $T\in \Q (G\times H)$ we have $|H|\cdot \tau_{G\times H} (T) = \tau_G (\iota(T))$. The result follows from Lemma \ref{lemma_trace_preserving_homo}.
 
\end{proof}

\begin{lemma}\label{lemma_additivity}
 Let $G$ be a countable discrete group. The set $\call C G$ is closed under addition. Furthermore, if $H$ is another countable discrete group, $a\in \call C (G)$, $b\in \call C(H)$ then $a+b \in \call C  (G\times H)$.
\end{lemma}
\begin{proof}
 The first claim follows from the fact that if $S\in M_k (\Q G)$ and $T\in M_l(\Q G)$  then $S\oplus T \in M_{k+l}(\Q G)$ has the property
$$
 \dimvn \ker (S\oplus T) = \dimvn \ker S + \dimvn \ker T.
$$ 

The second claims follows from taking $S\in M_k (\Q G)$ and $T\in M_l(\Q H)$ and observing that for $S\oplus T \in M_{k+l}(\Q (G\times H))$ we also have
$$
 \dimvn \ker (S\oplus T) = \dimvn \ker S + \dimvn \ker T.
$$ 
\end{proof}

\begin{lemma}\label{lemma_main_lemma} Let $(X,T_X)$ be a Turing dynamical system in which $X$ is a compact abelian group $\prod \Zmod{2}$, the action of $\Ga$ on $X$ is by continuous group automorphisms and the distinguished disjoint subsets $X_i$ of $X$ are cylinder sets.

Suppose furthermore that $(X,T_X)$ stops on any configuration, doesn't restart, and has disjoint accepting chains. Then $\mu(I) - \Om_1(X,T_X)$ is an $l^2$-Betti number arising from the group $\wh{X}\rtimes_{\wh{\rho}} \Ga$.
\end{lemma}

\begin{proof}
By Corollary \ref{main_cory}, $\mu(I) - \Om_1(X,T_X)$ is equal to $\dimvn \ker S$, where $S \in \call G(\rho)$ is expressed by a finite sum of elements $\ga_i\chi_i$, where $\ga_i \in \Ga$, and $\chi_i$ are products of characteristic functions of the sets $X_i$. By Lemma \ref{lemma_cylinder_sets}, $S$ is in the image of the Pontryagin map $P\otimes 1\colon  \Q(\wh X  \rtimes_{\wh \rho} \Ga) \to \call G (\rho)$. Let $\wh S$ be the preimage of $S$. By Proposition \ref{prop_pontryagin} we get 
$$
\dimvn \ker \wh S = \dimvn \ker S = \mu(I) - \Om_1(X,T_X).
$$
\end{proof}

\subsec{The lamplighter group}\label{subsec_proof_lamplighter}

\begin{theorem}
The set of $l^2$-Betti numbers arising from the group $(\Zmod{2}\wr \Z)^3$ contains 
$$
\frac1{64}- \frac18\sum_{k=1}^\infty \frac1{2^{k^2+4k+6}},
$$
\end{theorem}
\begin{proof}

By Proposition \ref{2oracle_reductions}, the system $(Y,T_Y)$ from Subsection \ref{section_tds-examples}\ref{subsec_second_example} fulfills conditions of Lemma \ref{lemma_main_lemma}. We conclude, by Corollary \ref{cory_2fundamental_value}, that
$$
	\frac1{64}- \frac18\sum_{k=1}^\infty \frac1{2^{k^2+4k+6}}
$$
is an $l^2$-Betti number arising from the group $\wh{Y}\rtimes_\wh{\rho} \De$, or more explicitly from
$$
	(\Zmod{2}^{\oplus \Z} \oplus \Zmod{2}^{\oplus \Z} \oplus \Zmod{2}^{\oplus \Z} \oplus \wh{S})
\rtimes_\wh{\rho} (\Z\oplus \Z\oplus \Z \oplus \Aut(S)).
$$
Note that the copies of $\Z$ act only on respective copies of $\Zmod{2}^{\oplus \Z}$, and that they act by the
shift. It follows that the above group is isomorphic to
$$
	(\Zmod{2}\wr \Z)^3 \times (\wh{S}\rtimes \Aut(S)),
$$
so the result follows from Lemma \ref{lemma_finite_group}.
\end{proof}

\subsec{Finitely generated groups}\label{subsec_proof_finitely_generated}

\begin{theorem*}
The set of $l^2$-Betti numbers arising from finitely generated groups is equal to the set of non-negative real numbers.
\end{theorem*}
\begin{proof}
By Proposition \ref{prop_tds_properties}, the system $(X^\Si,T_X^\Si)$ from Subsection \ref{section_tds-examples}.\ref{subsec_first_example} fulfills conditions of Lemma \ref{lemma_main_lemma}. We conclude, by Corollary \ref{cory_fundamental_value} that
$$
\frac1{64}- \frac2{8^3} \sum_{i\in \Si} \frac1{2^i}
$$ 
is an $l^2$-Betti number arising from the group $\wh X \rtimes_\wh{\rho^\Si} \Ga$, or more explicitly from
$$
	(\wh{M}^{\oplus \Z} \times \wh S)\,\, \rtimes_\wh{\rho}\,\, [(\Aut(M)\wr\Z) \ast \Zmod{2}]\times \Aut(S),
$$
so, by Lemma \ref{lemma_finite_group}, also from
\begin{equation}\label{eq_the_group}
   \wh{M}^{\oplus \Z}   \rtimes_\wh{\rho} ((\Aut(M)\wr\Z) \ast \Zmod{2}),
\end{equation}
which is easily seen to be finitely generated. This and the fact that $\Si$ was arbitrary implies that every
number between $\frac1{64}- \frac2{8^3}$ and $\frac1{64}$ is an $l^2$-Betti number arising from a finitely
generated group.

By additivity (Lemma \ref{lemma_additivity}), we see that there exists a natural number $N = \frac{8^3}{64\cdot 2}$ such that every number between $N -1$ and $N$ with is an $l^2$-Betti number arising from a finitely generated group. By additivity again it follows that every number bigger than $N$  is an $l^2$-Betti number arising from a finitely generated group. 

Let $r$ be an arbitrary positive real number. By the previous paragraph there exists a natural number $k$ such that $k \cdot r$ is an $l^2$-Betti number arising from a finitely presented group. The result follows from Lemma \ref{lemma_finite_group}. 
\end{proof}

\subsec{Finitely presented groups}\label{subsec_proof_finitely_presented}

\begin{theorem*}
The set of $l^2$-Betti numbers arising from finitely presented groups contains all numbers with computable binary expansions. 
\end{theorem*}

\begin{proof}
Recall that a presentation $\langle g_1, g_2, \ldots\,|\, r_1,r_2,\ldots\rangle$ is a {\bf recursive presentation} if there
exists an algorithm which lists all the elements of the set $\{r_1,r_2,\ldots\}$ in some order (perhaps with
repetitions). 

\begin{theorem*}[Higman's Embedding Theorem, \cite{higman:subgroups_of_finitely_presented_groups}] If a group has a recursive presentation then it can be
embedded into a finitely presented group.
\end{theorem*}

\begin{lemma*}\label{lemma_recursive_presentation} Let $\Si$ be a set of natural numbers and let $(X^\Si,T_X^\Si)$ be the Turing dynamical system from Subsection \ref{section_tds-examples}\ref{subsec_first_example}. If $\Si$ is computable then the group 
$$
   \wh{M}^{\oplus \Z}   \rtimes_\wh{\rho^\Si} ((\Aut(M)\wr\Z) \ast \Zmod{2})
$$
has a recursive presentation.
\end{lemma*}

\begin{proof} We use notation from Subsection \ref{section_tds-examples}\ref{subsec_first_example}. Recall that when a group $G = \langle g_1, g_2, \ldots \,|\, p_1, p_2, \ldots \rangle $ acts on a group $H = \langle h_1, h_2,\ldots \,|\, r_1,r_2, \ldots \rangle$ through an action $\al$ then the standard presentation of the
semi-direct product $H\rtimes_\al G$ is 
$$
\langle  h_2,h_2,\ldots ; g_1,g_2,\ldots \,|\, p_1,p_2,\ldots; r_1,r_2,\ldots; \,g_ih_jg_i^{-1} = \al(g_i)(h_j),\,\, i,j=1,2,\ldots \rangle.
$$

When we proceed to write this presentation in the case at hand, the only part which could possibly make it non-algorithmic is the action of $\Zmod{2}$. For $m\in \wh{\Zmod{2}^3}$ and $j\in \Z$, let $m_j$ denote the element $(\ldots,0,m,0,\ldots)\in
(\wh{\Zmod{2}^3})^{\oplus \Z}$, with $m$ on the $j$'th place. The relations which we have to write down are of the form:
$$
\begin{array}{l l l l}
	  B \cdot m_j \cdot B^{-1} &=& m_j, &\quad \text{ for } j\notin \Si	\\
	  B \cdot m_j \cdot B^{-1} &=& [\wh{\be'}(m)]_j &\quad \text{for } j\in \Si
\end{array}
$$
It is clear that if $\Si$ is computable, i.e. there is an algorithm which lists all elements of $\Si$ and there is an algorithm which lists all elements of the complement of $\Si$, then there exists also an algorithm which lists all of the above relations. 
\end{proof}

We saw in the previous subsection that the number
\begin{equation}\label{eq_the_number}
\frac1{64}- \frac2{8^3} \sum_{i\in \Si} \frac1{2^i}
\end{equation}
is an $l^2$-Betti arising from
$$
   \wh{M}^{\oplus \Z}   \rtimes_\wh{\rho^\Si} ((\Aut(M)\wr\Z) \ast \Zmod{2}).
$$
By Higman's Embedding Theorem and the previous lemma we see that for a computable $\Si$ the number  \eqref{eq_the_number} is an $l^2$-Betti number arising from a finitely presented group.

By additivity, there exists a natural number $N = \frac{8^3}{64\cdot 2}$ such that every number between $N -1$ and $N$ with a computable binary expansion is an $l^2$-Betti number arising from a finitely presented group. By additivity again, it follows that every number bigger than $N$ with a computable binary expansion is an $l^2$-Betti number arising from a finitely presented group. 

Let $r$ be an arbitrary positive number with a computable binary expansion. By the previous paragraph there exists a natural number $k$ such that $2^k \cdot r$ is an $l^2$-Betti number arising from a finitely presented group (since $2^k \cdot r$ has also computable binary expansion). The result follows from Lemma \ref{lemma_finite_group}. 
\end{proof}




\bibliographystyle{alpha}
\bibliography{bibliografia}
\end{document}